\pdfoutput=1
\RequirePackage{silence}
\documentclass[a4paper,svgnames]{amsart}
\usepackage[hmarginratio={1:1},vmarginratio={1:1},lmargin=60.0pt,tmargin=59.0pt]{geometry}

\usepackage{booktabs}
\allowdisplaybreaks

\usepackage[T1]{fontenc}
\usepackage{latexsym,exscale,amsfonts,amssymb,mathtools}
\usepackage{amsmath,amsthm,amsfonts,amssymb,amscd,textcomp,bbm}
\usepackage{mathrsfs,stackrel}
\usepackage{etoolbox}
\usepackage{tikz,tikz-cd}
\usetikzlibrary{matrix,arrows.meta,decorations.markings,shapes.multipart,decorations.pathreplacing,decorations.pathmorphing,shapes.geometric}
\usepackage{aliascnt}
\usepackage{ytableau}
\usepackage{xparse}
\usepackage{cite}

\usepackage{dynkin-diagrams}
\usepackage{nicematrix}
\usepackage{tikz-3dplot}
\usepackage{mathrsfs}
\DeclareMathAlphabet{\mathscrbf}{OMS}{mdugm}{b}{n}
\usepackage{float}
\usepackage{caption}
\usepackage{subcaption}
\usepackage[algoruled,lined]{algorithm2e}

\usepackage{array}
\newcolumntype{C}{>{$}c<{$}}

\usepackage{xcolor,colortbl}

\definecolor{lightgreen}{HTML}{CCFFCC}
\definecolor{lightblue}{HTML}{CCCCFF}

\definecolor{mygray}{gray}{0.6}
\definecolor{mygraydark}{gray}{0.4}
\definecolor{mygraylight}{gray}{0.85}
\definecolor{spinach}{RGB}{46,139,87}
\definecolor{tomato}{RGB}{255,99,71}
\definecolor{orchid}{RGB}{143,40,194}
\definecolor{neon}{RGB}{77,77,255}
\definecolor{lightneon}{RGB}{110,110,255}
\definecolor{pumpkin}{RGB}{224,180,80}
\definecolor{citron}{RGB}{190,180,90}

\definecolor{lava}{RGB}{207,16,32}
\definecolor{cream}{RGB}{255,253,208}
\definecolor{verdigris}{RGB}{67,179,174}
\definecolor{Black}{RGB}{0,0,0}
\definecolor{mydarkblue}{RGB}{10,10,170}
\definecolor{darkspinach}{RGB}{20,70,20}
\definecolor{darktomato}{RGB}{155,40,30}
\definecolor{darkorchid}{RGB}{50,10,100}
\definecolor{darklava}{RGB}{150,8,16}

\definecolor{zero}{RGB}{0,0,0}
\definecolor{one}{RGB}{255,0,0}
\definecolor{two}{RGB}{0,255,0}
\definecolor{three}{RGB}{0,0,255}

\usepackage{todonotes}

\usepackage{enumitem}
\setlist[enumerate]{itemsep=0.15cm,label=\emph{\upshape(\alph*)}}
\setlist[enumerate,2]{itemsep=0.15cm,label=\emph{\upshape(\roman*)}}
\setlist[enumerate,3]{itemsep=0.15cm,label=\emph{\upshape(\Alph*)}}

\let\emph\relax
\DeclareTextFontCommand{\emph}{\bfseries\em}

\renewcommand{\dots}{\text{...}}

\usepackage[all]{xy}
\usepackage{tikz}
\usetikzlibrary{cd}
\usetikzlibrary{decorations}
\usetikzlibrary{decorations.markings}
\usetikzlibrary{decorations.pathreplacing}
\usetikzlibrary{decorations.pathmorphing}
\usetikzlibrary{arrows.meta,shapes,positioning,matrix,calc}
\usetikzlibrary{shapes.callouts}
\tikzset{anchorbase/.style={baseline={([yshift=-0.5ex]current bounding box.center)}},
tinynodes/.style={font=\tiny,text height=0.25ex,text depth=0.05ex},
smallnodes/.style={font=\scriptsize,text height=0.75ex,text depth=0.15ex},
usual/.style={line width=2.0,color=black},
crossline/.style={preaction={draw=white,line width=2.75pt,-}},
}

\tikzstyle{startstop} = [ellipse, draw, fill=blue!20, text width=2cm, text centered, minimum height=1cm]
\tikzstyle{process} = [rectangle, draw, fill=orange!20, text width=2.5cm, text centered, minimum height=1cm]
\tikzstyle{arrow} = [thick, -{Stealth[length=3mm]}]

\newcommand{\R}{\mathbb{R}}

\newcommand{\K}{\mathbb{K}}

\def\NewTheorem#1{%
\newaliascnt{#1}{equation}%
\newtheorem{#1}[#1]{#1}%
\aliascntresetthe{#1}%
\expandafter\def\csname #1autorefname\endcsname{#1}%
}
\def\equationautorefname~#1\null{(#1)\null}

\numberwithin{equation}{subsection}

\NewTheorem{Proposition}
\NewTheorem{Theorem}
\NewTheorem{Corollary}
\AtEndEnvironment{Corollary}{\null\hfill$\square$}%
\NewTheorem{Lemma}
\NewTheorem{Conjecture}
\NewTheorem{Speculation}
\NewTheorem{Question}
\NewTheorem{Assumption}
\theoremstyle{definition}
\NewTheorem{Definition}
\AtEndEnvironment{Definition}{\null\hfill$\Diamond$}%
\NewTheorem{Classification Problem}
\AtEndEnvironment{Classification Problem}{\null\hfill$\Diamond$}%
\NewTheorem{Notation}
\AtEndEnvironment{Notation}{\null\hfill$\Diamond$}%
\NewTheorem{Example}
\AtEndEnvironment{Example}{\null\hfill$\Diamond$}%
\NewTheorem{Examples}
\AtEndEnvironment{Examples}{\vskip-10mm\null\hfill$\Diamond$}%

\theoremstyle{remark}
\NewTheorem{Remark}
\AtEndEnvironment{Remark}{\null\hfill$\Diamond$}%

\usepackage[hypertexnames=false]{hyperref}
\usepackage{bookmark}
\hypersetup{
pdftoolbar=true,
pdfmenubar=true,
pdffitwindow=false,
pdfstartview={FitH},
pdftitle={Representation Gap of the Motzkin Monoid},
pdfauthor={Katharina Arms},
pdfsubject={},
pdfcreator={Katharina Arms},
pdfproducer={Katharina Arms},
pdfkeywords={},
pdfnewwindow=true,
colorlinks=true,
linkcolor=mydarkblue,
citecolor=teal,
filecolor=magenta,
urlcolor=orchid,
linkbordercolor=lava,
citebordercolor=teal,
urlbordercolor=orchid,
linktocpage=true
}

\usepackage{listings}
\usepackage{xcolor}

\definecolor{codegray}{rgb}{0.95,0.95,0.95}
\definecolor{darkgreen}{rgb}{0,0.5,0}
\definecolor{darkblue}{rgb}{0,0,0.6}

\lstdefinelanguage{Mathematica}{
  morekeywords={Plot, Sum, Binomial, Mo},
  sensitive=true,
  morecomment=[s(*,*)],
  morestring=[b]"
}

\usepackage{courier}

\definecolor{codegray}{rgb}{0.95,0.95,0.95}
\definecolor{darkgreen}{rgb}{0,0.5,0}
\definecolor{darkblue}{rgb}{0,0,0.6}

\lstdefinelanguage{GAP}{
  morekeywords={
    function, local, return, if, then, fi, for, in, do, od,
    List, Length, Matrix, PrincipalFactor, ShallowCopy,
    DClasses, Rank, MotzkinMonoid
  },
  sensitive=true,
  morecomment=[l]\#,
}

\newcommand{\monoid}[1][S]{\mathcal{#1}}

\newcommand{\group}[1][G]{\mathcal{#1}}

\newcommand{\oneb}{\mathbbm{1}_{b}}
\newcommand{\onet}{\mathbbm{1}_{t}}
\newcommand{\onebt}{\mathbbm{1}_{bt}}

\newcommand{\lcell}{\mathcal{L}}
\newcommand{\rcell}{\mathcal{R}}
\newcommand{\jcell}{\mathcal{J}}
\newcommand{\hcell}{\mathcal{H}}

\newcommand{\lmod}[1][\lcell]{\Delta_{#1}}
\newcommand{\rmod}[1][\rcell]{{}_{#1}\Delta}

\def\makeautorefname#1#2{\csdef{#1autorefname}{#2}}
\makeautorefname{section}{Section}%
\makeautorefname{subsection}{Section}%
\makeautorefname{subsubsection}{Section}%

\begin{document}
\title[Representation Gap of the Motzkin Monoid]{Representation Gap of the Motzkin Monoid}
\author[K. Arms]{Katharina Arms}

\address{K.A.: School of Mathematics and Statistics, Carslaw Building (F07), University of Sydney, NSW 2006} 
\email{katharina.arms123@gmail.com}

\begin{abstract}
The linear decomposition attack reveals a vulnerability in encryption algorithms operating within groups or monoids with excessively small representations. The representation gap, defined as the size of the smallest non-trivial representation, therefore serves as a metric to assess the security of these algorithms. This paper will demonstrate that the diagrammatic Motzkin monoids exhibit a large representation gap, positioning them as promising candidates for robust encryption algorithms.
\end{abstract}


\addtocontents{toc}{\protect\setcounter{tocdepth}{1}}

\maketitle

\tableofcontents

\section{Introduction}


Many classical public-key protocols are instantiated in finite abelian groups. 
Diffie--Hellman key exchange, for example, is performed in a cyclic group 
$G=\langle g\rangle$ (typically $G=(\mathbb{F}_p)^\times$ or an elliptic-curve group), 
where parties choose secrets $a,b$, publish $g^a,g^b\in G$, and derive the shared 
secret $g^{ab}=g^{ba}$; commutativity underlies this equality (authentication is 
still required to prevent person-in-the-middle attacks); see \cite{Ko-algebraic-cryptography} for an overview.


Beyond commutative groups, active research explores cryptographic protocols based on 
non-commutative groups, monoids, semigroups, and even monoidal categories, motivated in 
part by the prospect of quantum computers threatening classical schemes; see 
\cite{MyShUs-group-cryptography,MyShUs-noncom-cryptography} for surveys. 
The guiding principle is to design systems whose security does not depend on 
hiding the underlying algebraic structure.


An early proposal of this type is Stickel's secret exchange 
\cite{St-new-key-exchange}, which can be viewed as a non-commutative analogue of 
Diffie--Hellman. Here one works in a non-commutative monoid $\mathcal{S}$ generated by 
elements $g,h\in \mathcal{S}$, together with commuting submonoids generated by powers 
of $g$ and of $h$. Parties $A$ and $B$ choose private integers $a,a'\in A$ and $b,b'\in B$, 
form $\alpha_p=g^a h^{a'}$ and $\beta_p=g^b h^{b'}$, and exchange these over a public channel. 
Each can then recover the common key;
\[
   g^b \alpha_p h^{b'} = g^a g^b h^{b'} h^{a'} = g^a \beta_p h^{a'},
\]
since powers of $g$ commute with each other, as do powers of $h$, while $g$ and $h$ do not. 
The security rests on the difficulty of recovering $(a,a')$ or $(b,b')$ from the 
public elements, $\alpha_p$ and $\beta_p$. This is a form of the decomposition problem, but its security \emph{heavily depends} on the used monoid.

A central motivation for studying non-traditional algebraic structures in 
cryptographic protocols is their susceptibility to the \emph{linear decomposition 
attack} \cite{MyRo-linear-attack}. This attack exploits finite-dimensional linear 
representations of the underlying algebraic object to break down the protocol. 

The \emph{representation gap} (repgap), the smallest dimension of a non-trivial 
finite-dimensional representation, provides a useful heuristic for gauging 
viability: if all non-trivial representations are of large dimension, then the 
protocol cannot be efficiently linearised and attacked by linear algebra. 
(Of course, the converse does not hold: a large repgap does not automatically 
make an algebraic structure cryptographically secure.)

In \cite{KST-monoidal-crypto}, Khovanov, Sitaraman and Tubbenhauer investigated 
diagrammatic monoids as candidates for Stickel’s and related non-commutative 
protocols, with the Temperley--Lieb monoid as a case study and the Motzkin monoid 
among the promising further examples. A key observation of 
\cite{KST-monoidal-crypto} is that finite groups often admit low-dimensional 
representations, enabling linear decomposition attacks and thereby forcing 
protocols to rely on increasingly large groups. This motivated the turn to 
monoids, particularly those arising from monoidal categories: by results of 
\cite{CouOsTub-growth}, and in line with Schur--Weyl duality phenomena, 
endomorphism monoids in such categories are expected to have exponential 
representation growth.

This paper focuses on one such example. Specifically, we establish a lower bound 
for the representation gap of the Motzkin monoid, showing it to be exponentially large. 

Our examination of the Motzkin monoid is useful and distinctive for the following reasons:
\begin{enumerate}
    \item By focusing on lower bounds rather than exact values, we obtain simple and 
    transparent asymptotic formulas that capture the essential growth behaviour, 
    whereas exact closed forms are often intractable. 

    \item The exponential growth in the representation gap of the Motzkin monoid is a promising initial indicator of its candidacy for resisting linear decomposition attacks in cryptographic applications. 

    \item The Motzkin monoid is defined diagrammatically, and our analysis likewise 
    employs diagrammatic methods. This makes our approach compatible with techniques 
    common in quantum topology and algebra, suggesting that the notion of 
    representation gaps for monoids may also have implications in these areas. 
\end{enumerate}

Methodologically, our work relies on the representation theory of monoids 
\cite{St-rep-monoid}, with Green's relations \cite{Gr-structure-semigroups} 
providing a natural cell structure. This framework both classifies and truncates 
elements of the Motzkin monoid and aligns with the formula counting Motzkin 
diagrams on $n$ nodes via a summation over $k$, the number of through strands, 
which index the cells. The use of Gram matrices to encode these cells then 
facilitates the study of ranks relative to idempotent elements, which is 
instrumental in bounding the representation gap. Altogether, this framework 
connects naturally with the theory of cellular algebras 
\cite{GrLe-cellular,Tu-sandwich}.

\subsection*{What is new in this paper}


The Motzkin monoid (or category) has long been studied 
(see, e.g., \cite{BeHa-motzkin}); however, its representation theory has 
received little attention, and in particular the dimensions of its simple 
representations remain unknown. In this paper we prove, for the first time, 
an exponential lower bound on the dimensions of its representations with 
base~3. This strengthens the earlier bound of base~2 obtained in 
\cite{KST-monoidal-crypto}, and, while broadly suggested by the main theorem of 
\cite{CouOsTub-growth}, is far from immediate. Our main result, stated in 
\autoref{T:Main}, thus determines the representation gap of the Motzkin monoid. 
Much of the preparatory material, such as parts of \autoref{S:Main}, is also new. 

\subsection*{Acknowledgments.}
I would first and foremost like to acknowledge my supervisor, Dani Tubbenhauer for their extensive advice and guidance in writing this paper. I would also like to acknowledge my mother, Kirsten, father, Chris and brother, Tommy for their support, and Amelie Skelton, Ellie Sage and Sarah Carstens for their comments.

This paper itself is a condensed version of the author's Honours thesis completed at the University of Sydney in 2025.

\section{Background}


Nothing in this section is new, and we refer to \cite{KST-monoidal-crypto} for details and references.

Throughout, let $\monoid$ be a finite monoid with group of units $\group$. For us,
the \emph{representation gap} $\operatorname{gap}_{\K}(\monoid)$ over a field $\K$ is defined as
\begin{gather*}
\operatorname{gap}_{\K}(\monoid) = \min \{ \dim_{\K}(L) | L\not\cong \onebt \hspace{0.1cm} \text{is a simple $\monoid$-representation}\}.
\end{gather*}
For our purposes, this definition is that same as in \cite{KST-monoidal-crypto}, as we will explain now.
Let $H^1(\group, \K)$ denote the standard group cohomology with values in $\K$. We call a monoid null-connected if every noninvertible element can be written as a nontrivial product of noninvertible elements. We call $\monoid$ right-connected if there is a unique equivalence class in $\monoid\setminus\group$ under $ab=a$ for $a,b\in\monoid\setminus\group$, and left-connected if the opposite monoid is right-connected. Finally, a monoid is \emph{well-connected} if it is a group or right-, left- and null-connected.

\begin{Theorem}
     \label{RepGap!}
     For a well-connected monoid $\monoid$ with $H^1(\group, \K) \cong 0$ the representation gap as defined above matches \cite[Deﬁnition 2A.6]{KST-monoidal-crypto}.
\end{Theorem}

\begin{proof}
    This is a direct reference of \cite[Theorem 2B.10]{KST-monoidal-crypto}.
\end{proof}

We will see in \autoref{S:Ext} that the conditions in \autoref{RepGap!} are satisfied for the monoid of our choice.

\begin{Definition}
    A \emph{faithful} representation is one on which every element of $\monoid$ acts uni\-quely, and the faithfulness of $\monoid$, denoted $\operatorname{faith}_{\K}(\monoid)$, is the minimum dimension of faithful $\monoid$ representations over $\K$.
\end{Definition}

\begin{Definition}
\label{ratios}
    We define the \emph{gap-ratio} and the \emph{faithful-ratio} of $\monoid$ as
    follows:
    \begin{equation*}
        \operatorname{gapr}_{\K}(\monoid) = \frac{\operatorname{gap}_{\K}(\monoid)}{\sqrt{ | \monoid |} }, \quad \operatorname{faithr}_{\K}(\monoid) = \frac{\operatorname{faith}_{\K}(\monoid)}{ \sqrt{| \monoid |} }.
    \end{equation*}
    (The normalisation by $\sqrt{|\monoid|}$ comes from the fact that simple $\monoid$-representations are at most of this size.)
\end{Definition}

We wish to find families of monoids, $ \{ \monoid_{n} \hspace{0.1cm} | \hspace{0.1cm} 0 \leq n \in \mathbb{Z} \} $ where these ratios do not approach 0 exponentially fast, for $n \rightarrow \infty$, and, at the same time, $\operatorname{gap}_{\K}(\monoid_n)$ grows exponentially.

\subsection{Cell Theory}
In order to define and study representations of monoids, we use \emph{Green's relations/cells}, which define equivalence classes that form group-like structures within the monoid. In order to construct the necessary equivalence classes, we begin by defining notions of left, right and two-sided order. For any monoid $\monoid$ and elements $a,b,c,d \in \monoid$,  we define these pre-orders as follows:
\[
(a \leq_{l} b) \Leftrightarrow \exists c :  b = ca,
\]
\[
(a \leq_{r} b) \Leftrightarrow \exists c :  b = ac,
\]
\[
(a \leq_{lr} b) \Leftrightarrow \exists c, d :  b = cad.
\]

\begin{Remark}
These orders align with those used in cellular algebras and in \cite{KST-monoidal-crypto}, but are reversed relative to the standard conventions in semigroup theory.
\end{Remark}

From these, we deduce notions of right, left and two-sided equivalence:
\[
(a \sim_{l} b) \Leftrightarrow (a \leq_{l} b \hspace{0.1cm} \text{and} \hspace{0.1cm} b \leq_{l} a),
\]
\[
(a \sim_{r} b) \Leftrightarrow (a \leq_{r} b \hspace{0.1cm} \text{and} \hspace{0.1cm} b \leq_{r} a),
\]
\[
(a \sim_{lr} b) \Leftrightarrow (a \leq_{lr} b \hspace{0.1cm} \text{and} \hspace{0.1cm} b \leq_{lr} a).
\]
We use $\lcell$, $\rcell$ and $\jcell$ to denote the above left, right and two-sided equivalence classes and call them the respective \emph{cells}. Furthermore, an H-cell is made up of elements in the intersection of a left cell $\lcell$ and right cell $\rcell$, denoted as $\hcell = \hcell(\lcell,\rcell) = \lcell \cap \rcell$. Refer to \cite{KST-monoidal-crypto} for details, and a visual representation. Later, we will classify monoids in terms of their J-cells. Note that the cells themselves are also ordered by the right, left and two-sided orders.

\begin{Lemma}
    There exist minimal and maximal elements of the $\leq_{lr}$ order, and they form the bottom and top J-cells $\jcell_{b}, \jcell_{t}$, respectively.
\end{Lemma}
\begin{proof}
    Refer to \cite{KST-monoidal-crypto} and \cite[Proposition 1.23]{St-rep-monoid}. 
\end{proof}

For our monoids we always have $|\hcell|=1$ for all H-cells, and we will assume this in the following to simplify the exposition. In this case, if $\hcell$ contains an idempotent, then we write $\hcell(e)$, which is the trivial group $1$.

\subsection{Classification of Simple Representations}
Cells are $\monoid$-representations, referred to as \emph{cell} or \emph{Schüt\-zenberger} representations and defined over $\K$.

\begin{Lemma}
    For any left cell $\lcell$, we have a left $\monoid$-representation which we denote by $\lmod[\lcell] = \K\lcell$, and is defined as follows:
    \begin{gather*}
        a \cdot l\in \Delta_{\lcell}=
        \begin{cases}
        al&\text{if }al\in\lcell,
        \\
        0&\text{else.}
        \end{cases}
    \end{gather*}
    Right representations $\rmod[\rcell]$ with respect to $\rcell$, and bi-representations with respect to $\jcell$, are defined similarly, and form the respective types of representations. 
\end{Lemma}
\begin{proof}
    Follows from definitions.
\end{proof}

The dimension of these are deduced directly from the definition to be:
\begin{equation*}
        \dim_{\K}(\lmod[\lcell]) = |\lcell|, \hspace{0.2cm} \dim_{\K}(\rmod[\rcell]) = |\rcell|.
\end{equation*}

The \emph{annihilator} of an $\monoid$-representation is the set of elements in $\monoid$ which send the representation to 0; denoted: $\operatorname{Ann}_{\monoid}(M) = \{ s \in \monoid \hspace{0.1cm}| \hspace{0.1cm} s \cdot M = 0 \}$.

\begin{Definition}
    An (or `the' momentarily) \emph{apex} of a cell-representation $M$ is the maximal J-cell $\jcell$ such that $\jcell$ contains no elements which annihilate $M$, and any J-cells $\jcell'$ which do must satisfy $\jcell' \leq_{lr} \jcell$. 
\end{Definition}

\begin{Lemma}
    The apex of any $\monoid$-representation $M$ exists and is unique.
\end{Lemma}

\begin{proof}
    This is well-known, see for example \cite{GaMaSt-irreps-semigroups}.
\end{proof}

We also recall that the maximal semisimple quotient of $M$, is referred to as its \emph{head}. We will denote this as $\operatorname{Hd}(M)$. We finally recall that the radical of $M$, $\operatorname{Rad}(M)$, is equal to the intersection of all of its maximal sub representations, and $\operatorname{Hd}(M) \cong M/\operatorname{Rad}(M)$, which is well defined for all $M$ over finite $\monoid$, up to isomorphism.

\begin{Proposition}(\emph{H-reduction}.)
\label{Clifford}
For any $\monoid$ and up to isomorphism, there is one simple $\monoid$-representation for every apex $\jcell$, and a J-cell $\jcell$ is an apex if and only if it contains some $\hcell(e)$.

Moreover, if $\jcell$ is an apex and $\lcell\subset\jcell$ is any left cell containing some $\hcell(e)$, then the associated simple $\monoid$-representation is $L=\operatorname{Hd}(\lmod[\lcell])$.
\end{Proposition}

\begin{proof}
    This can be proven via interpretation of \cite[Theorem 5.5]{St-rep-monoid}, where we define isomorphism classes by J-cells.
\end{proof}

There is also a right cell version of \autoref{Clifford}, but we will not need it.

\subsubsection{Cells and Semisimple Representation Gaps}

We now consider the dimension of representations.

\begin{Definition}
    \label{ssdim}
    Let $\jcell$ be an apex with associated simple $L$.
    We now define the \emph{semisimple dimension} of $L$ to be:
    \[
    \operatorname{ssdim} (L) =\operatorname{dim}_{\K}(\lmod[\lcell]) = |\lcell|.
    \]
    The \emph{semisimple representation gap} is the smallest $m$ for which there is a non-trivial semisimple $\monoid$-representation with $\operatorname{ssdim}(L) = m$.
    The \emph{semisimple-gap-ratio} is:
    \[\operatorname{ssgapr}(\monoid) = \frac{\operatorname{ssgap}(\monoid)}{\sqrt{| \monoid |}},\]
    where we again normalize by $\sqrt{| \monoid |}$. This is characteristic free and even independent of $\K$, so we use no $\K$ as a subscript.
\end{Definition}

\subsubsection{Cells and Gram Matrices}

Recall that $\lcell$ and $\rcell$ are free right and left sets with respective sets of representatives; $\{ l_{1}, ..., l_{R} \}$ and $\{ r_{1}, ..., r_{L} \}$, where R and L are the number of right and left cells in $\jcell$. We denote the \emph{Gram matrix} by $P^{\jcell} = (P_{i,j}^{\jcell})_{i,j}$, a matrix over $\K \hcell(e)$, with entries: 
\[
P_{i,j}^{\jcell} = 
\begin{cases}
1 \quad \text{if}  \hspace{0.1cm} r_{i}l_{j}=q1, \\
0 \quad else.
\end{cases}
\]

\begin{Proposition}
    \label{rank}
    For cells within a fixed J-cell, assuming that $P^{\jcell}$ is square and symmetric, we have:
        \[
        \dim_{\K}(L_{\jcell}) = \operatorname{rank}(P^{\jcell}).
        \]
    Where we denote the associated simple $\monoid$-represen\-tation by $L_{\jcell}$.
\end{Proposition}

\begin{proof}
    Refer to \cite{KST-monoidal-crypto} for details of this proof.
\end{proof}

Call $\jcell$ idempotent if it contains some $\hcell(e)$.

\begin{Theorem}
    \label{dimtruncate}
    Holding the above assumptions, for a submonoid $\mathcal{R} \subset \monoid$, if $\jcell$ restricts to an idempotent J-cell of $\mathcal{R}$, then for the associated simple $\mathcal{R}$ and $\monoid$ representations:
    \[
    \dim_{\K}(L_{\jcell}^{\monoid}) \geq \dim_{\K}(L_{\jcell}^{\mathcal{R}})
    \]
\end{Theorem}

\begin{proof}
    Since $\mathcal{R}$ is a submonoid, this follows from the fact that the size of a monoid is always larger than or equivalent to that of any of its submonoids.
\end{proof}



\subsection{Cell Submonoids}\label{Submonoid}

We see later in \autoref{subquotient} that the growth in representation gap is observed mostly in a small number of cells, motivating us to construct a submonoid restricted to J-cells within a certain range.

\begin{Lemma}
    \label{lsubquotient}
    For any J cell $\jcell$, the set $\monoid_{\geq\jcell}=\{s\in\monoid|s\geq_{lr}\jcell\}$ is an ideal in $\monoid$.
\end{Lemma}
\begin{proof}
    Easy and omitted.
\end{proof}

Using these, we construct our \emph{cell sub-quotients}:

\begin{Definition}
    The $\mathcal{K}$-submonoid, for a given J-cell, $\mathcal{K}$, of some $\monoid$ is a monoid.
    \[
    \monoid^{\mathcal{K}} = \monoid_{\geq\mathcal{K}}\cup\{1\},
    \]
    where we attach a formal unit (if necessary).
    This is a monoid by \autoref{lsubquotient}.
\end{Definition}

We assume for the rest of this section that $\monoid$ is regular (meaning every J-cell is an apex), and can define the J-cells of its sub-quotient as follows:

\begin{Lemma}
    The J-cells of $\monoid^{\mathcal{K}}$ can be expressed as such:
    \[
    \{ \jcell_{b} \} \cup \{ \mathcal{M} \hspace{0.1cm} | \hspace{0.1cm} \mathcal{M}  \hspace{0.1cm} \text{is a J-cell of $\monoid$ for} \hspace{0.1cm} \mathcal{M} \geq_{lr} \mathcal{K} \}.
    \]
\end{Lemma}
\begin{proof}
See \cite[Lemma 3F.7]{KST-monoidal-crypto}.
\end{proof}

Representations of $\jcell_{b}$ and $\jcell_{t}$ in $\monoid^{\mathcal{K}}$ are uniquely given by $\oneb$ and $\onet$ respectively. For representations of $\mathcal{M} \not\cong \jcell_{b}, \jcell_{t}$, we note that $\monoid^{\mathcal{K}}$ representations can be \emph{inflated} to $\monoid$ representations, giving the following correspondence:

\begin{Proposition}
    \label{truncate}
    If $\mathcal{M} \not\in \{\jcell_{b}, \jcell_{t} \}$ is an apex of  $\monoid^{\mathcal{K}}$ and $\monoid$, then inflation makes the simple $\monoid^{\mathcal{K}}$-representa\-tion of apex $\mathcal{M}$ into a $\monoid$-representation of apex $\mathcal{M}$ (in particular, they have the same dimension).
\end{Proposition}
\begin{proof}
    The first equation is a direct result of \autoref{Clifford}, and the dimensions of representations will not change since $\monoid$ is regular, and thus both represent analogous sets.
\end{proof}

Thus, we have the following bound for simple $\monoid$-representations:

\begin{Theorem}
    \label{lower bound on truncation}
    Given a J-cell, $\mathcal{K}$ of $\monoid$, 
    \[
    \operatorname{gap}_{\K}(\monoid^{\mathcal{K}}) \geq \operatorname{gap}_{\K}(\monoid),\quad
    \operatorname{ssgap}(\monoid^{\mathcal{K}}) \geq \operatorname{ssgap}(\monoid).
    \]
\end{Theorem}
\begin{proof}
    This follows directly from \autoref{truncate} above.
\end{proof}

\section{The Motzkin Monoid and Its Properties}\label{S:Main}

We will now prove a lower bound for the representation gap of the Motzkin monoid, based on the mathematical groundwork from the previous section.

\subsection{Defining the Motzkin Monoid}

The Motzkin monoid, $\mathcal{M}o_n$, has submonoids $\mathcal{TL}_n$ (Temperley--Lieb) and $p\mathcal{R}o_n$ (planar rook), and can be defined in a similar way. It has been dealt with in numerous papers, including \cite{DoEaGr-motzkin-pbrauer} and \cite{BeHa-motzkin}. It also has the same number of elements as the Motzkin numbers, which count the number of ways one can draw crossingless lines between nodes on a circle (not necessarily touching every node), see for example \cite{Callan-CatalanMotzkinRiordan}. 
These numbers are easily seen to grow 
exponentially, with exponential factor 9.
Despite this, extensive research on the potential for the Motzkin monoid to be used in cryptography remains to be explored, so we begin in this direction by proving its promising representation gap (of square root size of the Motzkin numbers), which motivates further exploration. Our results were used in \cite{SteTub-rigid-monoids}.

\vspace{0.2cm}

The following description of the Motzkin monoid is formulated with reference to the construction of the Temperley--Lieb monoid in \cite{KST-monoidal-crypto}, working over an arbitrary field $\K$. This is a set-theoretical version of the classical Motzkin algebra as in, for example, \cite{BeHa-motzkin}. Below, we describe the space of endomorphisms, which we use to define diagrammatic elements of the Motzkin monoid. 

\vspace{0.2cm}

We will begin by considering the ($\K$ linear) Motzkin category, which we denote $\mathcal{M}o^{lin}(\delta)$, and consists of nodes along the strip $\R \times [0,1]$, with $n$ nodes existing along the bottom i.e. when $t = 0$, and $m$ nodes along the top, when $t = 1$ (where $t \in [0,1]$). The set of morphisms between these nodes are given in the Motzkin category by all crossingless partitions of the bottom and top points with at most two parts. That is, nodes may be paired to nodes in the other row, to nodes in their own row, or be disconnected. Throughout the paper, we use the term \emph{through strands} to refer to strands which connect nodes in the bottom and top row, i.e. matchings between a node at $t=0$ and at $t=1$. These matchings are represented by diagrams, the following being examples:
\begin{gather*}
    \begin{tikzpicture}[anchorbase]
    \draw[usual] (0.25,0) to[out=90,in=180] (0.5,0.25) to[out=0,in=90] (0.75,0);
    \draw[usual] (0.75,1) to[out=270,in=180] (1,0.75) to[out=0,in=270] (1.25,1);
    \draw[usual] (0,1) to (0,0);
    \draw[usual] (0.25,1) to (1.25,0);
    \end{tikzpicture}
    \; \hspace{1cm} \;
    \begin{tikzpicture}[anchorbase]
    \draw[usual] (0,1) to[out=270,in=180] (0.25,0.75) to[out=0,in=270] (0.5,1);
    \draw[usual] (0,0) to (0,0.25);
    \draw[usual] (0.5,0) to (0.5,0.25);
    \end{tikzpicture}
    \; \hspace{1cm} \;
    \begin{tikzpicture}[anchorbase]
    \draw[usual] (0,0) to[out=90,in=180] (0.5,0.5) to[out=0,in=90] (1,0);
    \draw[usual] (0.5,1.5) to[out=270,in=180] (1.25,1) to[out=0,in=270] (2,1.5);
    \draw[usual] (1,1.5) to[out=270,in=180] (1.25,1.25) to[out=0,in=270] (1.5,1.5);
    \draw[usual] (0,1.5) to (1.5,0);
    \draw[usual] (2,0) to (2,0.5);
    \draw[usual] (0.5,0) to (0.5,0.25);
    \end{tikzpicture}
    \; \hspace{1cm} \;
    \begin{tikzpicture}[anchorbase]
    \draw[usual] (1.5,0) to (1.5,1);
    \draw[usual] (0,0) to (0,0.25);
    \draw[usual] (1,1) to (1,0.75);
    \draw[usual] (0,1) to (0.5,0);
    \draw[usual] (0.5,1) to (1,0);
    \end{tikzpicture}
\end{gather*}
(We draw partitions of size one as ``ending lines'' to make them visually more prominent.)
Here we understand crossingless to mean that it is possible to represent any element of $\mathcal{M}o^{lin}(\delta)$ with a planar graph. Also, each node may only be attached to one line. For example, the following are not elements of $\mathcal{M}o^{lin}(\delta)$:
\begin{gather*}
\text{Nonexamples: }
    \begin{tikzpicture}[anchorbase]
    \draw[usual] (0.5,0) to[out=90,in=180] (0.75,0.5) to[out=0,in=90] (1,0);
    \draw[usual] (0,1) to (0,0);
    \draw[usual] (0.5,1) to (0,0);
    \draw[usual] (0.5,1) to (0.5,0);
    \draw[usual] (1,1) to (1,0);
    \draw[usual] (1.5,1) to (1,0);
    \draw[usual] (1.5,0.5) to (1.5,0);
    \end{tikzpicture}
    \; \hspace{1cm} \;
    \begin{tikzpicture}[anchorbase]
    \draw[usual] (0.5,1) to[out=270,in=180] (0.75,0.5) to[out=0,in=270] (1,1);
    \draw[usual] (1,1) to[out=270,in=180] (1.5,0.5) to[out=0,in=270] (2,1);
    \draw[usual] (0,0) to (0,0.5);
    \draw[usual] (0.5,0) to (0.5,1);
    \draw[usual] (0,1) to (1,0);
    \draw[usual] (1.5,0) to (1.5,1);
    \draw[usual] (2,0) to (1.5,1);
    \end{tikzpicture}
\end{gather*}
Composition in
$\mathcal{M}o^{lin}(\delta)$ is furthermore defined via the operation of `gluing', whereby objects are stacked on top of each other to create another element (conventions are specified in the picture below). It is possible for composition to result in an element which has components which are not connected to nodes in the top or bottom row. These lines may appear as `squiggles' (intervals topologically) or closed circles. We evaluate them to some chosen $\delta\in\K$ as in the following example:
\begin{gather*}
    \begin{tikzpicture}[anchorbase]
    \draw[usual] (0,0) to[out=90,in=180] (0.25,0.25) to[out=0,in=90] (0.5,0);
    \draw[usual] (0.5,1) to[out=270,in=180] (0.75,0.75) to[out=0,in=270] (1,1);
    \draw[usual] (1,0) to[out=90,in=180] (1.25,0.25) to[out=0,in=90] (1.5,0);
    \draw[usual] (0,1) to (0,0.75);
    \draw[usual] (1.5,1) to (1.5,0.75);
    \end{tikzpicture}
    \; \circ \;
    \begin{tikzpicture}[anchorbase]
    \draw[usual] (0,1) to[out=270,in=180] (0.25,0.75) to[out=0,in=270] (0.5,1);
    \draw[usual] (1,1) to[out=270,in=180] (1.25,0.75) to[out=0,in=270] (1.5,1);
    \draw[usual] (0,0) to[out=90,in=180] (0.5,0.5) to[out=0,in=90] (1,0);
    \draw[usual] (0.5,0) to (0.5,0.25);
    \draw[usual] (1.5,0) to (1.5,0.25);
    \end{tikzpicture}
    \; \hspace{0.2cm} = \hspace{0.2cm} \;
    \begin{tikzpicture}[anchorbase]
    \draw[usual] (0,1) to[out=90,in=180] (0.25,1.25) to[out=0,in=90] (0.5,1);
    \draw[usual] (0.5,2) to[out=270,in=180] (0.75,1.75) to[out=0,in=270] (1,2);
    \draw[usual] (1,1) to[out=90,in=180] (1.25,1.25) to[out=0,in=90] (1.5,1);
    \draw[usual] (0,2) to (0,1.75);
    \draw[usual] (1.5,2) to (1.5,1.75);
    \draw[usual] (0,1) to[out=270,in=180] (0.25,0.75) to[out=0,in=270] (0.5,1);
    \draw[usual] (1,1) to[out=270,in=180] (1.25,0.75) to[out=0,in=270] (1.5,1);
    \draw[usual] (0,0) to[out=90,in=180] (0.5,0.5) to[out=0,in=90] (1,0);
    \draw[usual] (0.5,0) to (0.5,0.25);
    \draw[usual] (1.5,0) to (1.5,0.25);
    \end{tikzpicture}
    \; \hspace{0.2cm} = \hspace{0.2cm} \delta^2 \cdot \;
    \begin{tikzpicture}[anchorbase]
    \draw[usual] (0.5,1) to[out=270,in=180] (0.75,0.75) to[out=0,in=270] (1,1);
    \draw[usual] (0,1) to (0,0.75);
    \draw[usual] (1.5,1) to (1.5,0.75);
    \draw[usual] (0,0) to[out=90,in=180] (0.5,0.5) to[out=0,in=90] (1,0);
    \draw[usual] (0.5,0) to (0.5,0.25);
    \draw[usual] (1.5,0) to (1.5,0.25);
    \end{tikzpicture}
\end{gather*}
By definition, we evaluate to $\delta^n$ for $n$ disconnected objects. The monoidal structure on the Motzkin category is defined as in Temperley--Lieb, namely that the operator $\otimes$ represents horizontal adding of diagrams; $m \otimes n = m+n$. The example below can easily be generalised to understand this convention:
\begin{gather*}
    \begin{tikzpicture}[anchorbase]
    \draw[usual] (0,1) to (0.5,0);
    \draw[usual] (1,1) to (1,0);
    \draw[usual] (0,0) to (0,0.25);
    \draw[usual] (0.5,1) to (0.5,0.75);
    \end{tikzpicture}
    \; \otimes \;
    \begin{tikzpicture}[anchorbase]
    \draw[usual] (0.5,1) to[out=270,in=180] (0.75,0.75) to[out=0,in=270] (1,1);
    \draw[usual] (0,0) to[out=90,in=180] (0.25,0.25) to[out=0,in=90] (0.5,0);
    \draw[usual] (0,1) to (1,0);
    \draw[usual] (1.5,0) to (1.5,1);
    \end{tikzpicture}
    \; \hspace{0.2cm} = \hspace{0.2cm} \;
    \begin{tikzpicture}[anchorbase]
    \draw[usual] (0,1) to (0.5,0);
    \draw[usual] (1,1) to (1,0);
    \draw[usual] (0,0) to (0,0.25);
    \draw[usual] (0.5,1) to (0.5,0.75);
    \end{tikzpicture}
    \; \;
    \begin{tikzpicture}[anchorbase]
    \draw[usual] (0.5,1) to[out=270,in=180] (0.75,0.75) to[out=0,in=270] (1,1);
    \draw[usual] (0,0) to[out=90,in=180] (0.25,0.25) to[out=0,in=90] (0.5,0);
    \draw[usual] (0,1) to (1,0);
    \draw[usual] (1.5,0) to (1.5,1);
    \end{tikzpicture}
\end{gather*}
The Motzkin algebra $\mathcal{M}o^{lin}_n(\delta)$ on n strands is the endomorphism algebra of n nodes.

We now define generators of the Motzkin algebra, as constructed in \cite{PosHatLy-MotzkinMonoid}. The algebra $\mathcal{M}o^{lin}_n(\delta)$ has a presentation with three generators, $k_i$, $r_i$ and $l_i$, as follows, each of these enumerating by $i$, denoting a pattern occurring at the $i^{th}$ node of a Motzkin element:
\begin{gather*}
    \underbrace{\begin{tikzpicture}[anchorbase]
    \draw[usual] (0.5,0) to (0.5,0.25);
    \draw[usual] (0,1) to (0,0.75);
    \draw[usual] (0,0) to (0.5,1);
    \end{tikzpicture}}_{l_i}
    \;\hspace{0.2cm}\;
    \underbrace{\begin{tikzpicture}[anchorbase]
    \draw[usual] (0,1) to[out=270,in=180] (0.25,0.75) to[out=0,in=270] (0.5,1);
    \draw[usual] (0,0) to[out=90,in=180] (0.25,0.25) to[out=0,in=90] (0.5,0);
    \end{tikzpicture}}_{t_i}
    \;\hspace{0.2cm}\;
    \underbrace{\begin{tikzpicture}[anchorbase]
    \draw[usual] (0.5,0) to (0.5,0.25);
    \draw[usual] (0,0) to (0.5,1);
    \draw[usual] (0,1) to (0,0.75);
    \end{tikzpicture}}_{r_i}
\end{gather*}

\begin{Lemma}
    The Motzkin matchings between a row of $m$ nodes and one of $n$ nodes which form a $\K$-linear basis for the homomorphisms from $m$ to $n$ objects. The size of this basis is given by the Motzkin number $\sum_{k = 0}^n\frac{1}{k+1}\binom{2n}{2k}\binom{2k}k{}$ where $k$ represents the number of possible through strands ranging from 0 to $n$.
\end{Lemma}

\begin{proof}
    This is non-trivial, but well-known, see e.g. \cite{BeHa-motzkin}.
\end{proof}

Note here that unlike in Temperley--Lieb, we do not require $\frac{m+n}{2}$ to be an integer, as an uneven number of nodes can exist in a Motzkin diagram, with at least one node disconnected.

\begin{Lemma}
    \label{antiinvolution}
    The Motzkin monoid category has an anti-involution, given by reflection over the central horizontal axis. We denote the anti-involution of an element $s \in \mathcal{M}o$ as $s^*$. An example in $\mathcal{M}o_5$ is pictured below:
\end{Lemma}

\begin{gather*}
    \left( \hspace{0.1cm}
    \begin{tikzpicture}[anchorbase]
    \draw[usual] (0,0) to[out=90,in=180] (0.25,0.25) to[out=0,in=90] (0.5,0);
    \draw[usual] (0.5,1) to[out=270,in=180] (0.75,0.75) to[out=0,in=270] (1,1);
    \draw[usual] (0,1) to (1,0);
    \draw[usual] (2,1) to (1.5,0);
    \draw[usual] (2,0) to (2,0.25);
    \draw[usual] (1.5,1) to (1.5,0.75);
    \end{tikzpicture}
    \hspace{0.1cm} \right)^{\ast}
    \; = \;
    \begin{tikzpicture}[anchorbase]
    \draw[usual] (0,1) to[out=270,in=180] (0.25,0.75) to[out=0,in=270] (0.5,1);
    \draw[usual] (0.5,0) to[out=90,in=180] (0.75,0.25) to[out=0,in=90] (1,0);
    \draw[usual] (0,0) to (1,1);
    \draw[usual] (2,0) to (1.5,1);
    \draw[usual] (2,1) to (2,0.75);
    \draw[usual] (1.5,0) to (1.5,0.25);
    \end{tikzpicture}
\end{gather*}

\begin{proof}
    Clear by definition; swap top and bottom presentations and this must also be a morphism of our category, and doing it twice gives back the original diagram.
\end{proof}

When considering cryptographic applications of the Motzkin monoid, the existence of $\delta$ to signify an unmatched circle or squiggle is not useful. As such, we will consider all $\delta$'s to be trivial, i.e. $\delta=1$, and disregard them in the composition as in example below. 

\begin{gather*}
    \begin{tikzpicture}[anchorbase]
    \draw[usual] (0.5,0) to[out=90,in=180] (0.75,0.25) to[out=0,in=90] (1,0);
    \draw[usual] (1,1) to[out=270,in=180] (1.25,0.75) to[out=0,in=270] (1.5,1);
    \draw[usual] (0,1) to (0,0);
    \draw[usual] (0.5,1) to (1.5,0);
    \end{tikzpicture}
    \; \circ \;
    \begin{tikzpicture}[anchorbase]
    \draw[usual] (0.5,1) to[out=270,in=180] (0.75,0.75) to[out=0,in=270] (1,1);
    \draw[usual] (1.5,0) to (1.5,0.25);
    \draw[usual] (0,0) to (0,0.25);
    \draw[usual] (0,1) to (0.5,0);
    \draw[usual] (1.5,1) to (1,0);
    \end{tikzpicture}
    \; \hspace{0.2cm} = \hspace{0.2cm} \;
    \begin{tikzpicture}[anchorbase]
    \draw[usual] (0.5,1) to[out=90,in=180] (0.75,1.25) to[out=0,in=90] (1,1);
    \draw[usual] (1,2) to[out=270,in=180] (1.25,1.75) to[out=0,in=270] (1.5,2);
    \draw[usual] (0,2) to (0,1);
    \draw[usual] (0.5,2) to (1.5,1);
    \draw[usual] (0.5,1) to[out=270,in=180] (0.75,0.75) to[out=0,in=270] (1,1);
    \draw[usual] (1.5,0) to (1.5,0.25);
    \draw[usual] (0,0) to (0,0.25);
    \draw[usual] (0,1) to (0.5,0);
    \draw[usual] (1.5,1) to (1,0);
    \end{tikzpicture}
    \; \hspace{0.2cm} = \hspace{0.2cm} \delta \cdot \;
    \begin{tikzpicture}[anchorbase]
    \draw[usual] (1,1) to[out=270,in=180] (1.25,0.75) to[out=0,in=270] (1.5,1);
    \draw[usual] (1.5,0) to (1.5,0.25);
    \draw[usual] (0,0) to (0,0.25);
    \draw[usual] (0,1) to (0.5,0);
    \draw[usual] (0.5,1) to (1,0);
    \end{tikzpicture}
    \; \hspace{0.2cm} = \hspace{0.2cm} \;
    \begin{tikzpicture}[anchorbase]
    \draw[usual] (1,1) to[out=270,in=180] (1.25,0.75) to[out=0,in=270] (1.5,1);
    \draw[usual] (1.5,0) to (1.5,0.25);
    \draw[usual] (0,0) to (0,0.25);
    \draw[usual] (0,1) to (0.5,0);
    \draw[usual] (0.5,1) to (1,0);
    \end{tikzpicture}
\end{gather*}

Formally, we specialize $\delta$ to $1$. The resulting category $\mathcal{M}o^{lin}(1)$ can then be treated completely set-theoretically, i.e. without any linear structure. Denote this category by $\mathcal{M}o$.

\vspace{0.2cm}

We only consider diagrams for which $m=n$, and consequently define the Motzkin monoid as the set of endomorphisms $\mathcal{M}o_n = \operatorname{End}_{\mathcal{M}o}(n)$.

\subsection{Representations Using Cell Theory}

The cells of the Motzkin monoid can be identified in a similar way to how \cite{KST-monoidal-crypto} does it for Temperley--Lieb. Here, each $\jcell_k$ consists of diagrams with the same number of through strands $k$ in the matchings, and each $\lcell$ and $\rcell$ within a given $\jcell$ contain elements with the same presentation on the bottom and top rows of nodes respectively. The picture below shows the cells of $\mathcal{M}o_3$, with idempotent elements depicted in red:
\begin{gather*}\label{Eq:TLCells}
\begin{gathered}
\xy
(0,0)*{\begin{gathered}
\begin{tabular}{C!{\color{blue}\vrule width 0.1mm}C!{\color{blue}\vrule width 0.1mm}C!{\color{blue}\vrule width 0.1mm}C}
\arrayrulecolor{blue}
\begin{tikzpicture}[anchorbase]
\draw[red, ultra thick] (0,0) to (0,0.25);
\draw[red, ultra thick] (0.5,0) to (0.5,0.25);
\draw[red, ultra thick] (1,0) to (1,0.25);
\draw[red, ultra thick] (0,0.75) to (0,0.5);
\draw[red, ultra thick] (0.5,0.75) to[out=270,in=180] (0.75,0.5) to[out=0,in=270] (1,0.75);
\end{tikzpicture} &
\begin{tikzpicture}[anchorbase]
\draw[red, ultra thick] (0,0) to (0,0.25);
\draw[red, ultra thick] (0.5,0.75) to[out=270,in=180] (0.75,0.5) to[out=0,in=270] (1,0.75);
\draw[red, ultra thick] (0.5,0) to[out=90,in=180] (0.75,0.25) to[out=0,in=90] (1,0);
\draw[red, ultra thick] (0,0.75) to (0,0.5);
\end{tikzpicture} &
\begin{tikzpicture}[anchorbase]
\draw[red, ultra thick] (1,0) to (1,0.25);
\draw[red, ultra thick] (0.5,0.75) to[out=270,in=180] (0.75,0.5) to[out=0,in=270] (1,0.75);
\draw[red, ultra thick] (0,0) to[out=90,in=180] (0.25,0.25) to[out=0,in=90] (0.5,0);
\draw[red, ultra thick] (0,0.75) to (0,0.5);
\end{tikzpicture} &
\begin{tikzpicture}[anchorbase]
\draw[red, ultra thick] (0.5,0.75) to[out=270,in=180] (0.75,0.5) to[out=0,in=270] (1,0.75);
\draw[red, ultra thick] (0,0) to[out=90,in=180] (0.25,0.25) to[out=0,in=90] (1,0);
\draw[red, ultra thick] (0,0.75) to (0,0.5);
\draw[red, ultra thick] (0.5,0) to (0.5,0.125);
\end{tikzpicture}
\\
\hline
\begin{tikzpicture}[anchorbase]
\draw[red, ultra thick] (0,0) to (0,0.25);
\draw[red, ultra thick] (0.5,0) to (0.5,0.25);
\draw[red, ultra thick] (1,0) to (1,0.25);
\draw[red, ultra thick] (1,0.75) to (1,0.5);
\draw[red, ultra thick] (0,0.75) to[out=270,in=180] (0.25,0.5) to[out=0,in=270] (0.5,0.75);
\end{tikzpicture} &
\begin{tikzpicture}[anchorbase]
\draw[red, ultra thick] (0,0) to (0,0.25);
\draw[red, ultra thick] (0.5,0) to[out=90,in=180] (0.75,0.25) to[out=0,in=90] (1,0);
\draw[red, ultra thick] (1,0.75) to (1,0.5);
\draw[red, ultra thick] (0,0.75) to[out=270,in=180] (0.25,0.5) to[out=0,in=270] (0.5,0.75);
\end{tikzpicture} &
\begin{tikzpicture}[anchorbase]
\draw[red, ultra thick] (1,0) to (1,0.25);
\draw[red, ultra thick] (0,0) to[out=90,in=180] (0.25,0.25) to[out=0,in=90] (0.5,0);
\draw[red, ultra thick] (1,0.75) to (1,0.5);
\draw[red, ultra thick] (0,0.75) to[out=270,in=180] (0.25,0.5) to[out=0,in=270] (0.5,0.75);
\end{tikzpicture} &
\begin{tikzpicture}[anchorbase]
\draw[red, ultra thick] (0,0) to[out=90,in=180] (0.25,0.25) to[out=0,in=90] (1,0);
\draw[red, ultra thick] (0.5,0) to (0.5,0.125);
\draw[red, ultra thick] (1,0.75) to (1,0.5);
\draw[red, ultra thick] (0,0.75) to[out=270,in=180] (0.25,0.5) to[out=0,in=270] (0.5,0.75);
\end{tikzpicture}
\\
\hline
\begin{tikzpicture}[anchorbase]
\draw[red, ultra thick] (0,0) to (0,0.25);
\draw[red, ultra thick] (0.5,0) to (0.5,0.25);
\draw[red, ultra thick] (1,0) to (1,0.25);
\draw[red, ultra thick] (0.5,0.75) to (0.5,0.625);
\draw[red, ultra thick] (0,0.75) to[out=270,in=180] (0.5,0.5) to[out=0,in=270] (1,0.75);
\end{tikzpicture} &
\begin{tikzpicture}[anchorbase]
\draw[red, ultra thick] (0,0) to (0,0.25);
\draw[red, ultra thick] (0.5,0) to[out=90,in=180] (0.75,0.25) to[out=0,in=90] (1,0);
\draw[red, ultra thick] (0.5,0.75) to (0.5,0.625);
\draw[red, ultra thick] (0,0.75) to[out=270,in=180] (0.5,0.5) to[out=0,in=270] (1,0.75);
\end{tikzpicture} &
\begin{tikzpicture}[anchorbase]
\draw[red, ultra thick] (1,0) to (1,0.25);
\draw[red, ultra thick] (0,0) to[out=90,in=180] (0.25,0.25) to[out=0,in=90] (0.5,0);
\draw[red, ultra thick] (0.5,0.75) to (0.5,0.625);
\draw[red, ultra thick] (0,0.75) to[out=270,in=180] (0.5,0.5) to[out=0,in=270] (1,0.75);
\end{tikzpicture} &
\begin{tikzpicture}[anchorbase]
\draw[red, ultra thick] (0,0) to[out=90,in=180] (0.25,0.25) to[out=0,in=90] (1,0);
\draw[red, ultra thick] (0.5,0) to (0.5,0.125);
\draw[red, ultra thick] (0.5,0.75) to (0.5,0.625);
\draw[red, ultra thick] (0,0.75) to[out=270,in=180] (0.5,0.5) to[out=0,in=270] (1,0.75);
\end{tikzpicture}
\\
\hline
\begin{tikzpicture}[anchorbase]
\draw[red, ultra thick] (0,0) to (0,0.25);
\draw[red, ultra thick] (0.5,0) to (0.5,0.25);
\draw[red, ultra thick] (1,0) to (1,0.25);
\draw[red, ultra thick] (0,0.75) to (0,0.5);
\draw[red, ultra thick] (0.5,0.75) to (0.5,0.5);
\draw[red, ultra thick] (1,0.75) to (1,0.5);
\end{tikzpicture} &
\begin{tikzpicture}[anchorbase]
\draw[red, ultra thick] (0,0) to (0,0.25);
\draw[red, ultra thick] (0.5,0) to[out=90,in=180] (0.75,0.25) to[out=0,in=90] (1,0);
\draw[red, ultra thick] (0,0.75) to (0,0.5);
\draw[red, ultra thick] (0.5,0.75) to (0.5,0.5);
\draw[red, ultra thick] (1,0.75) to (1,0.5);
\end{tikzpicture} &
\begin{tikzpicture}[anchorbase]
\draw[red, ultra thick] (1,0) to (1,0.25);
\draw[red, ultra thick] (0,0) to[out=90,in=180] (0.25,0.25) to[out=0,in=90] (0.5,0);
\draw[red, ultra thick] (0,0.75) to (0,0.5);
\draw[red, ultra thick] (0.5,0.75) to (0.5,0.5);
\draw[red, ultra thick] (1,0.75) to (1,0.5);
\end{tikzpicture} &
\begin{tikzpicture}[anchorbase]
\draw[red, ultra thick] (0,0) to[out=90,in=180] (0.25,0.25) to[out=0,in=90] (1,0);
\draw[red, ultra thick] (0.5,0) to (0.5,0.125);
\draw[red, ultra thick] (0,0.75) to (0,0.5);
\draw[red, ultra thick] (0.5,0.75) to (0.5,0.5);
\draw[red, ultra thick] (1,0.75) to (1,0.5);
\end{tikzpicture}
\end{tabular}
\\[3pt]
\begin{tabular}{C!{\color{blue}\vrule width 0.1mm}C!{\color{blue}\vrule width 0.1mm}C!{\color{blue}\vrule width 0.1mm}C!{\color{blue}\vrule width 0.1mm}C}
\arrayrulecolor{blue}
\begin{tikzpicture}[anchorbase]
\draw[ red,ultra thick] (0,0.75) to (0,0);
\draw[ red,ultra thick] (0.5,0.75) to[out=270,in=180] (0.75,0.5) to[out=0,in=270] (1,0.75);
\draw[ red,ultra thick] (0.5,0) to[out=90,in=180] (0.75,0.25) to[out=0,in=90] (1,0);
\end{tikzpicture} & 
\begin{tikzpicture}[anchorbase]
\draw[ red,ultra thick] (0,0.75) to (0,0);
\draw[ red,ultra thick] (0.5,0.75) to[out=270,in=180] (0.75,0.5) to[out=0,in=270] (1,0.75);
\draw[ red,ultra thick] (0.5,0) to (0.5,0.25);
\draw[ red,ultra thick] (1,0) to (1,0.25);
\end{tikzpicture} &
\begin{tikzpicture}[anchorbase]
\draw[usual] (0,0.75) to (0.5,0);
\draw[usual] (0.5,0.75) to[out=270,in=180] (0.75,0.5) to[out=0,in=270] (1,0.75);
\draw[usual] (0,0) to (0,0.25);
\draw[usual] (1,0) to (1,0.25);
\end{tikzpicture} &
\begin{tikzpicture}[anchorbase]
\draw[ red,ultra thick] (0,0.75) to (1,0);
\draw[ red,ultra thick] (0.5,0.75) to[out=270,in=180] (0.75,0.5) to[out=0,in=270] (1,0.75);
\draw[ red,ultra thick] (0,0) to (0,0.2);
\draw[ red,ultra thick] (0.5,0) to (0.5,0.2);
\end{tikzpicture} &
\begin{tikzpicture}[anchorbase]
\draw[usual] (0,0.75) to (1,0);
\draw[usual] (0.5,0.75) to[out=270,in=180] (0.75,0.5) to[out=0,in=270] (1,0.75);
\draw[usual] (0,0) to[out=90,in=180] (0.25,0.25) to[out=0,in=90] (0.5,0);
\end{tikzpicture} 
\\
\hline
\begin{tikzpicture}[anchorbase]
\draw[ red,ultra thick] (0,0.75) to (0,0);
\draw[ red,ultra thick] (0.5,0.75) to (0.5,0.5);
\draw[ red,ultra thick] (1,0.75) to (1,0.5);
\draw[ red,ultra thick] (0.5,0) to[out=90,in=180] (0.75,0.25) to[out=0,in=90] (1,0);
\end{tikzpicture} & 
\begin{tikzpicture}[anchorbase]
\draw[ red,ultra thick] (0,0.75) to (0,0);
\draw[ red,ultra thick] (0.5,0.75) to (0.5,0.5);
\draw[ red,ultra thick] (1,0.75) to (1,0.5);
\draw[ red,ultra thick] (0.5,0) to (0.5,0.25);
\draw[ red,ultra thick] (1,0) to (1,0.25);
\end{tikzpicture} &
\begin{tikzpicture}[anchorbase]
\draw[usual] (0,0.75) to (0.5,0);
\draw[usual] (0.5,0.75) to (0.5,0.5);
\draw[usual] (1,0.75) to (1,0.5);
\draw[usual] (0,0) to (0,0.25);
\draw[usual] (1,0) to (1,0.25);
\end{tikzpicture} &
\begin{tikzpicture}[anchorbase]
\draw[usual] (0,0.75) to (1,0);
\draw[usual] (0.5,0.75) to (0.5,0.55);
\draw[usual] (1,0.75) to (1,0.55);
\draw[usual] (0,0) to (0,0.2);
\draw[usual] (0.5,0) to (0.5,0.2);
\end{tikzpicture} &
\begin{tikzpicture}[anchorbase]
\draw[usual] (0,0.75) to (1,0);
\draw[usual] (0.5,0.75) to (0.5,0.55);
\draw[usual] (1,0.75) to (1,0.55);
\draw[usual] (0,0) to[out=90,in=180] (0.25,0.25) to[out=0,in=90] (0.5,0);
\end{tikzpicture} 
\\
\hline
\begin{tikzpicture}[anchorbase]
\draw[usual] (0.5,0.75) to (0,0);
\draw[usual] (0,0.75) to (0,0.5);
\draw[usual] (1,0.75) to (1,0.5);
\draw[usual] (0.5,0) to[out=90,in=180] (0.75,0.25) to[out=0,in=90] (1,0);
\end{tikzpicture} & 
\begin{tikzpicture}[anchorbase]
\draw[usual] (0.5,0.75) to (0,0);
\draw[usual] (0,0.75) to (0,0.5);
\draw[usual] (1,0.75) to (1,0.5);
\draw[usual] (0.5,0) to (0.5,0.25);
\draw[usual] (1,0) to (1,0.25);
\end{tikzpicture} &
\begin{tikzpicture}[anchorbase]
\draw[ red,ultra thick] (0.5,0.75) to (0.5,0);
\draw[ red,ultra thick] (0,0.75) to (0,0.5);
\draw[ red,ultra thick] (1,0.75) to (1,0.5);
\draw[ red,ultra thick] (0,0) to (0,0.25);
\draw[ red,ultra thick] (1,0) to (1,0.25);
\end{tikzpicture} &
\begin{tikzpicture}[anchorbase]
\draw[usual] (0.5,0.75) to (1,0);
\draw[usual] (0,0.75) to (0,0.55);
\draw[usual] (1,0.75) to (1,0.55);
\draw[usual] (0,0) to (0,0.2);
\draw[usual] (0.5,0) to (0.5,0.2);
\end{tikzpicture} &
\begin{tikzpicture}[anchorbase]
\draw[usual] (0.5,0.75) to (1,0);
\draw[usual] (0,0.75) to (0,0.55);
\draw[usual] (1,0.75) to (1,0.55);
\draw[usual] (0,0) to[out=90,in=180] (0.25,0.25) to[out=0,in=90] (0.5,0);
\end{tikzpicture}
\\
\hline
\begin{tikzpicture}[anchorbase]
\draw[ red,ultra thick] (1,0.75) to (0,0);
\draw[ red,ultra thick] (0,0.75) to (0,0.55);
\draw[ red,ultra thick] (0.5,0.75) to (0.5,0.55);
\draw[ red,ultra thick] (0.5,0) to[out=90,in=180] (0.75,0.25) to[out=0,in=90] (1,0);
\end{tikzpicture} & 
\begin{tikzpicture}[anchorbase]
\draw[usual] (1,0.75) to (0,0);
\draw[usual] (0.5,0.75) to (0.5,0.55);
\draw[usual] (0,0.75) to (0,0.55);
\draw[usual] (0.5,0) to (0.5,0.2);
\draw[usual] (1,0) to (1,0.2);
\end{tikzpicture} &
\begin{tikzpicture}[anchorbase]
\draw[usual] (1,0.75) to (0.5,0);
\draw[usual] (0,0.75) to (0,0.5);
\draw[usual] (0.5,0.75) to (0.5,0.5);
\draw[usual] (0,0) to (0,0.25);
\draw[usual] (1,0) to (1,0.25);
\end{tikzpicture} &
\begin{tikzpicture}[anchorbase]
\draw[ red,ultra thick] (1,0.75) to (1,0);
\draw[ red,ultra thick] (0,0.75) to (0,0.55);
\draw[ red,ultra thick] (0.5,0.75) to (0.5,0.55);
\draw[ red,ultra thick] (0,0) to (0,0.2);
\draw[ red,ultra thick] (0.5,0) to (0.5,0.2);
\end{tikzpicture} &
\begin{tikzpicture}[anchorbase]
\draw[ red,ultra thick] (1,0.75) to (1,0);
\draw[ red,ultra thick] (0,0.75) to (0,0.55);
\draw[ red,ultra thick] (0.5,0.75) to (0.5,0.55);
\draw[ red,ultra thick] (0,0) to[out=90,in=180] (0.25,0.25) to[out=0,in=90] (0.5,0);
\end{tikzpicture}
\\
\hline
\begin{tikzpicture}[anchorbase]
\draw[usual] (1,0.75) to (0,0);
\draw[usual] (0,0.75) to[out=270,in=180] (0.25,0.5) to[out=0,in=270] (0.5,0.75);
\draw[usual] (0.5,0) to[out=90,in=180] (0.75,0.25) to[out=0,in=90] (1,0);
\end{tikzpicture} & 
\begin{tikzpicture}[anchorbase]
\draw[usual] (1,0.75) to (0,0);
\draw[usual] (0,0.75) to[out=270,in=180] (0.25,0.5) to[out=0,in=270] (0.5,0.75);
\draw[usual] (0.5,0) to (0.5,0.25);
\draw[usual] (1,0) to (1,0.25);
\end{tikzpicture} &
\begin{tikzpicture}[anchorbase]
\draw[usual] (1,0.75) to (0.5,0);
\draw[usual] (0,0.75) to[out=270,in=180] (0.25,0.5) to[out=0,in=270] (0.5,0.75);
\draw[usual] (0,0) to (0,0.25);
\draw[usual] (1,0) to (1,0.25);
\end{tikzpicture} &
\begin{tikzpicture}[anchorbase]
\draw[ red,ultra thick] (1,0.75) to (1,0);
\draw[ red,ultra thick] (0,0.75) to[out=270,in=180] (0.25,0.5) to[out=0,in=270] (0.5,0.75);
\draw[ red,ultra thick] (0,0) to (0,0.2);
\draw[ red,ultra thick] (0.5,0) to (0.5,0.2);
\end{tikzpicture} &
\begin{tikzpicture}[anchorbase]
\draw[ red,ultra thick] (1,0.75) to (1,0);
\draw[ red,ultra thick] (0,0.75) to[out=270,in=180] (0.25,0.5) to[out=0,in=270] (0.5,0.75);
\draw[ red,ultra thick] (0,0) to[out=90,in=180] (0.25,0.25) to[out=0,in=90] (0.5,0);
\end{tikzpicture}
\end{tabular}
\\[3pt]
\begin{tabular}{C!{\color{blue}\vrule width 0.1mm}C!{\color{blue}\vrule width 0.1mm}C}
\arrayrulecolor{blue}
\begin{tikzpicture}[anchorbase]
\draw[red, ultra thick] (0,0) to (0,0.75);
\draw[red, ultra thick] (0.5,0) to (0.5,0.75);
\draw[red, ultra thick] (1,0.75) to (1,0.5);
\draw[red, ultra thick] (1,0) to (1,0.25);
\end{tikzpicture} &
\begin{tikzpicture}[anchorbase]
\draw[usual] (0,0) to (0,0.75);
\draw[usual] (1,0) to (0.5,0.75);
\draw[usual] (1,0.75) to (1,0.5);
\draw[usual] (0.5,0) to (0.5,0.25);
\end{tikzpicture} &
\begin{tikzpicture}[anchorbase]
\draw[usual] (0.5,0) to (0,0.75);
\draw[usual] (1,0) to (1,0.75);
\draw[usual] (1,0.75) to (1,0.5);
\draw[usual] (0,0) to (0,0.25);
\draw[usual] (0.5,0.75) to (0.5,0.5);
\end{tikzpicture}
\\
\hline
\begin{tikzpicture}[anchorbase]
\draw[usual] (0,0) to (0,0.75);
\draw[usual] (0.5,0) to (1,0.75);
\draw[usual] (0.5,0.75) to (0.5,0.5);
\draw[usual] (1,0) to (1,0.25);
\end{tikzpicture} &
\begin{tikzpicture}[anchorbase]
\draw[red, ultra thick] (0,0) to (0,0.75);
\draw[red, ultra thick] (1,0) to (1,0.75);
\draw[red, ultra thick] (0.5,0.75) to (0.5,0.5);
\draw[red, ultra thick] (0.5,0) to (0.5,0.25);
\end{tikzpicture} &
\begin{tikzpicture}[anchorbase]
\draw[usual] (0.5,0) to (0,0.75);
\draw[usual] (1,0) to (1,0.75);
\draw[usual] (0.5,0.75) to (0.5,0.5);
\draw[usual] (0,0) to (0,0.25);
\end{tikzpicture}
\\
\hline
\begin{tikzpicture}[anchorbase]
\draw[usual] (0,0) to (0.5,0.75);
\draw[usual] (0.5,0) to (1,0.75);
\draw[usual] (0,0.75) to (0,0.5);
\draw[usual] (1,0) to (1,0.25);
\end{tikzpicture} &
\begin{tikzpicture}[anchorbase]
\draw[usual] (0,0) to (0.5,0.75);
\draw[usual] (1,0) to (1,0.75);
\draw[usual] (0,0.75) to (0,0.5);
\draw[usual] (0.5,0) to (0.5,0.25);
\end{tikzpicture} &
\begin{tikzpicture}[anchorbase]
\draw[red, ultra thick] (0.5,0) to (0.5,0.75);
\draw[red, ultra thick] (1,0) to (1,0.75);
\draw[red, ultra thick] (0,0.75) to (0,0.5);
\draw[red, ultra thick] (0,0) to (0,0.25);
\end{tikzpicture}
\end{tabular}
\\[3pt]
\begin{tabular}{C}
\arrayrulecolor{blue}
\begin{tikzpicture}[anchorbase]
\draw[red, ultra thick] (0,0) to (0,0.5);
\draw[red, ultra thick] (0.5,0) to (0.5,0.5);
\draw[red, ultra thick] (1,0) to (1,0.5);
\end{tikzpicture}
\end{tabular}
\end{gathered}};
(-45,40)*{\jcell_{0}};
(-45,-0.5)*{\jcell_{1}};
(-45,-36.5)*{\jcell_{2}};
(-45,-53.5)*{\jcell_{3}};
\endxy
\quad.
\end{gathered}
\end{gather*}

\begin{Lemma}\label{L:Factor}
    There exists a unique factorisation for each element $a \in \operatorname{Hom}(m,n)$, namely, $a = b \circ i_k \circ c$ for $b \in \operatorname{Hom}(n,k)$, $i_k = id_k \in \operatorname{Hom}(k,k)$ and $c \in \operatorname{Hom}(k,m)$.
\end{Lemma}

\begin{proof}
    The following diagram exemplifies one such factorisation $a$:
    \begin{gather*}\label{Eq:TLFactor}
a=
\begin{tikzpicture}[anchorbase]
\draw[usual] (0.5,0) to[out=90,in=180] (1,0.5) to[out=0,in=90] (1.5,0);
\draw[usual] (0,1) to[out=270,in=180] (0.25,0.75) to[out=0,in=270] (0.5,1);
\draw[usual] (0,0) to (1,1);
\draw[usual] (1,0) to (1,0.25);
\draw[usual] (1.5,1) to (1.5,0.75);
\end{tikzpicture}
\;=\;
\underbrace{\begin{tikzpicture}[anchorbase]
\draw[usual] (0,1) to[out=270,in=180] (0.25,0.75) to[out=0,in=270] (0.5,1);
\draw[usual] (1,0) to (1,1);
\draw[usual] (1.5,1) to (1.5,0.75);
\end{tikzpicture}}_{b}
\;\circ\;
\underbrace{\begin{tikzpicture}[anchorbase]
\draw[usual] (0,0) to (0,1);
\end{tikzpicture}}_{\mathbf{id}_{1}}
\;\circ\;
\underbrace{\begin{tikzpicture}[anchorbase]
\draw[usual] (0.5,0) to[out=90,in=180] (1,0.5) to[out=0,in=90] (1.5,0);
\draw[usual] (0,0) to (0,1);
\draw[usual] (1,0) to (1,0.25);
\end{tikzpicture}}_{c}
\;,
\end{gather*}
    In this example, $a \in \mathcal{M}o_4$ has been factorised as described. This illustrates, without loss of generality, that factorisation would clearly be possible for any $s \in \mathcal{M}o$ generally.
\end{proof}

\begin{Proposition}
    \label{motzkin formula}
    For $\mathcal{M}o_n$, the map that assigns a diagram to its number of through strands gives 
    an order preserving one-to-one correspondence between J-cells and $\{n,\dots,0\}$ (in this reversed order). Furthermore, in any J-cell with $k$ through strands; $\jcell_k$, the number of left and right cells for the Motzkin monoid are equal, and are given by:
    \[
    |\lcell | = | \rcell | = \sum_{t=0}^{n}\frac{k + 1}{k+t+1}\binom{n}{k+2t}\binom{k+2t}{t}.
    \]
    We also have,
    \[
    | \jcell | =  | \lcell | \cdot | \rcell | = | \lcell |^2.
    \]
    
    All J-cells are idempotent, and for all idempotent H-cells, $\hcell(e) \cong 1$, where $| \hcell | = 1$.
\end{Proposition}

\begin{proof}
    The statement about the J-cells follows directly from \autoref{L:Factor}, along with the observation that the number of through strands cannot increase during composition.

    For the L-cells, the first formula is taken directly from \cite{KST-monoidal-crypto}. Since the left cells are an anti-involution of the right cells, we know that our cells will be square, and the second formula follows.

    Now, all J-cells are idempotent since vertically symmetric diagrams such as the following;
    \begin{gather*}
    \begin{tikzpicture}[anchorbase]
    \draw[usual] (0.5,0) to (0.5,1);
    \draw[usual] (0,0) to (0,1);
    \draw[usual] (1,1) to[out=270,in=180] (1.5,0.7) to[out=0,in=270] (2,1);
    \draw[usual] (1,0) to[out=90,in=180] (1.5,0.3) to[out=0,in=90] (2,0);
    \draw[usual] (1.5,1) to (1.5,0.8);
    \draw[usual] (1.5,0) to (1.5,0.2);
    \draw[usual] (2.5,1) to (2.5,0);
    \end{tikzpicture}
    \; \hspace{1cm} \;
    \begin{tikzpicture}[anchorbase]
    \draw[usual] (0.5,0) to (0.5,0.3);
    \draw[usual] (0.5,1) to (0.5,0.7);
    \draw[usual] (0,0) to (0,1);
    \draw[usual] (1,1) to (1,0);
    \draw[usual] (1.5,0) to (1.5,1);
    \end{tikzpicture}
    \; \hspace{1cm} \;
    \begin{tikzpicture}[anchorbase]
    \draw[usual] (0,1) to[out=270,in=180] (0.25,0.7) to[out=0,in=270] (0.5,1);
    \draw[usual] (0,0) to[out=90,in=180] (0.25,0.3) to[out=0,in=90] (0.5,0);
    \draw[usual] (1,1) to (1,0);
    \draw[usual] (1.5,1) to (1.5,0.7);
    \draw[usual] (1.5,0) to (1.5,0.3);
    \draw[usual] (2,0) to (2,1);
    \end{tikzpicture}
    \end{gather*}
    
    are idempotents, as one easily sees.
    
    The final statement about $\hcell$ is true since they contain exactly one element, namely that determined specifically by a unique bottom and top presentation, so that when this is idempotent, it forms a trivial group.
\end{proof}

\begin{Proposition}\label{P:Simples}
    The set of apexes of the Motzkin monoid can be indexed by the set 
    $\{n,\dots,0\}$ in a one-to-one correspondence where $n$ denotes the number of through strands of elements in a $\jcell$. These apexes correspond to a unique simple representation up to isomorphism.
\end{Proposition}

\begin{proof}
    This follows directly from H-reduction, \autoref{Clifford} and \autoref{motzkin formula}.
\end{proof}

We denote the simple representations from \autoref{P:Simples} by $L_k$ for $k$ the number of through strands and the apex of $L_k$.

\subsection{Sub-quotients: Truncating the Motzkin Monoid.}
\label{subquotient}
As in the Temperley--Lieb case, the Motzkin monoid only has a sufficiently high representation gap in cells with the optimal number of through strands. As shown in \autoref{strandsdiagram} below, high dimensional representations clearly only occur for low numbers of strands.

\begin{figure}[hbt!]
    \centering
    \includegraphics[scale = 0.45]{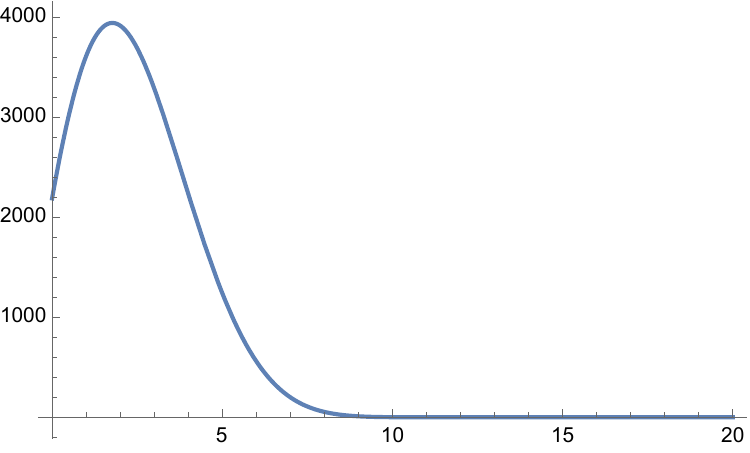} 
    \includegraphics[scale = 0.45]{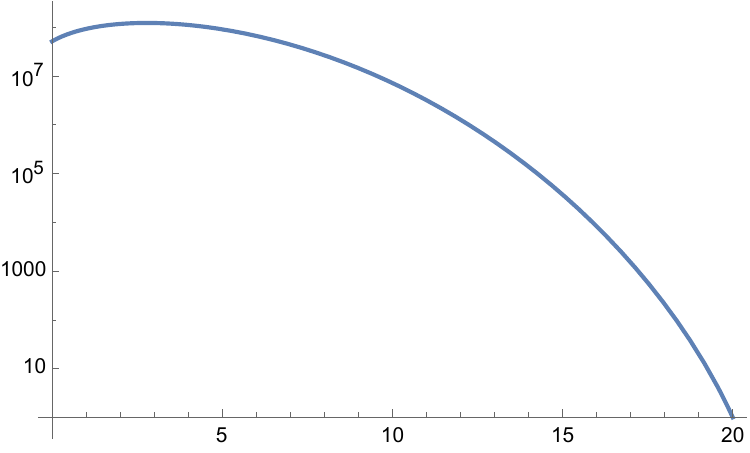}
    \caption{Number of diagrams existing for $k \in\{0,\dots,20\}$ over $n=20$ nodes. (Thus illustrating the semisimple dimensions.) The second plot uses a logarithmic scale on the y-axis.}
    \label{strandsdiagram}
\end{figure}

In many other cases, for example where $n \approx k$, all associated simple representations are small and we would like to omit them in order to prevent easy attacks on these cells. As such, we truncate the Motzkin monoid, using the process as described in \autoref{Submonoid}. We begin by proving that the Motzkin monoid meets the necessary conditions so that the cell structure remains intact after truncation.

\begin{Definition}
    Define the $k$-th truncated Motzkin monoid by:
    \[
    \mathcal{M}o_n^{\leq k} = (\mathcal{M}o_n)_{\geq \jcell_k}
    \]
    (We could also truncate from both sides, but the binomial-type peak in \autoref{strandsdiagram} suggests that truncation on one side is good enough.)
\end{Definition}

The Motzkin monoid is regular by \autoref{motzkin formula}, so this truncation is well-behaved with respect to the cell structure. The k-th truncated Motzkin monoid forms a monoid with at most $k$ through strands, as per the definition of cells and equivalence classes outlined in \autoref{Submonoid}. Referring to \autoref{strandsdiagram}, we can see that it would be optimal to allow for roughly $\sqrt{n}$ through strands for a monoid with $2n$ ($n$ in each row) nodes in order to maximize the representation gap, since the peak clearly occurs before this truncation. 

\begin{Remark}
In fact, by differentiating the formula for the size of left cells in \autoref{motzkin formula} one can show that the peak is always at $\leq\sqrt{n}$ through strands, so our truncation below will keep the peak.
\end{Remark}

\subsection{Trivial Extensions}\label{S:Ext}

In this section, we examine the behaviour of trivial extensions using notions of left, right and well-connectedness, as well as considering the actions of trivial elements. This is important in being able to prove our desired formula for the representation gap.

\begin{Lemma}\label{P:Split}
    The Motzkin monoid $M$ is well-connected, and has $H^1(\monoid, \K) \cong 0.$
\end{Lemma}

\begin{proof}
That the Motzkin monoid is null-connected follows (almost) directly from \autoref{L:Factor}.

    We next show that the monoid is left-connected. To be left-connected means that if we choose an element within the monoid, it can be obtained from any other element within the monoid via left-action by some $s \in \monoid$. We will prove this by considering the generators of our monoid, and showing them to be left-connected.

    We already have that the generators $t_i$ are connected via the concept of flip pairs as defined in \cite{KST-monoidal-crypto} for $n \geq 5$ and $k \geq 3$, so it remains to show that the other generators are also connected. 

    We can easily show that any $r_i$ and $l_i$ are left-connected. To do so, take arbitrary $r_m$ and $l_n$. Now consider $d = r_mr_ml_nl_n$, the diagram with disconnected nodes at $m,n,m+1,n+1$. We have $dr_m = d = dl_n$ showing that $r_i \approx_l l_i, \forall i$, and $l_i \approx_l r_i, \forall i$. 

    It now remains to show that some $t_i$ is left-connected to some $r_j$, and then left-connectedness of all generators will follow, since it is an equivalence relation. Namely, if all three generators on the $i$th node are left equivalent, and $t_i$ is left equivalent to all positions $t_n$, then $r_i$ and $l_i$ are equivalent to all $t_n$, and thus $r_n$ and $l_n$. Equivalence of $t_i$ and $r_j$ (for arbitrary non-close $i,j$) is also easy since $t_ir_j = r_jt_i$ We see this in the following diagrams:
    \begin{gather*}
    \begin{tikzpicture}[anchorbase]
    \draw[usual] (0,0) to (0,1);
    \draw[usual] (0.25,0) to (0.25,0.5);
    \draw[usual] (0.5,0) to (0.5,0.5);
    \draw[usual] (0.75,0) to (0.75,1);
    \draw[usual] (0.25,1) node[above]{$i$} to[out=270,in=180] (0.375,0.85) to[out=0,in=270] (0.5,1);
    \draw[usual] (0.25,0.5) to[out=90,in=180] (0.375,0.65) to[out=0,in=90] (0.5,0.5);
    \end{tikzpicture}
    \;\hspace{0.2cm} \cdots \hspace{0.2cm}\;
    \begin{tikzpicture}[anchorbase]
    \draw[usual] (0,0) to (0,1);
    \draw[usual] (0.25,0) to (0.5,0.5);
    \draw[usual] (0.5,0) to (0.5,0.15);
    \draw[usual] (0.75,0) to (0.75,1);
    \draw[usual] (0.5,0.5) to (0.5,1);
    \draw[usual] (0.25,0.35) to (0.25,1) node[above]{$j$};
    \end{tikzpicture}
    \;\hspace{0.2cm} = \hspace{0.2cm}\;
    \begin{tikzpicture}[anchorbase]
    \draw[usual] (0,0) to (0,0.5);
    \draw[usual] (0.75,0) to (0.75,0.5);
    \draw[usual] (0.25,0.5) node[above]{$i$} to[out=270,in=180] (0.375,0.35) to[out=0,in=270] (0.5,0.5);
    \draw[usual] (0.25,0) to[out=90,in=180] (0.375,0.15) to[out=0,in=90] (0.5,0);
    \end{tikzpicture}
    \;\hspace{0.2cm} \cdots \hspace{0.2cm}\;
    \begin{tikzpicture}[anchorbase]
    \draw[usual] (0,0) to (0,0.5);
    \draw[usual] (0.25,0) to (0.5,0.5);
    \draw[usual] (0.5,0) to (0.5,0.15);
    \draw[usual] (0.75,0) to (0.75,0.5);
    \draw[usual] (0.25,0.35) to (0.25,0.5) node[above]{$j$};
    \end{tikzpicture}
    \;\hspace{0.2cm} = \hspace{0.2cm}\;
    \begin{tikzpicture}[anchorbase]
    \draw[usual] (0,0) to (0,1);
    \draw[usual] (0.25,0.5) to (0.25,1) node[above]{$i$};
    \draw[usual] (0.5,0.5) to (0.5,1);
    \draw[usual] (0.75,0) to (0.75,1);
    \draw[usual] (0.25,0.5) to[out=270,in=180] (0.375,0.35) to[out=0,in=270] (0.5,0.5);
    \draw[usual] (0.25,0) to[out=90,in=180] (0.375,0.15) to[out=0,in=90] (0.5,0);
    \end{tikzpicture}
    \;\hspace{0.2cm} \cdots \hspace{0.2cm}\;
    \begin{tikzpicture}[anchorbase]
    \draw[usual] (0,0) to (0,1);
    \draw[usual] (0.25,0) to (0.25,0.5);
    \draw[usual] (0.5,0) to (0.5,0.65);
    \draw[usual] (0.75,0) to (0.75,1);
    \draw[usual] (0.5,1) to (0.25,0.5);
    \draw[usual] (0.25,0.85) to (0.25,1) node[above]{$j$};
    \end{tikzpicture}
\end{gather*}
    This proves that these are left-connected for $k\geq 2,n \geq 4$. Hence, all generators are left-connected, and since each element can be expressed as a product of generators, we have proven that the Motzkin monoid is left-connected. Due to the anti-involution defined in \autoref{antiinvolution}, left-connected\-ness implies right-connectedness, and we have hence proven that the Motz\-kin monoid is well-connected.

    In order to see that $H^1(\monoid, \K) \cong 0$, we refer to the proof of \cite[Lemma 4D.15]{KST-monoidal-crypto}. Note that an equivalent argument applies here, since all elements in the Motzkin monoid have a finite period.
\end{proof}

We can now state our definition of the representation gap to be the minimum simple dimension as a theorem for this section.
\begin{Theorem}
    \label{simple=gap}
    The representation gap of the Motzkin monoid, as defined in \cite[Deﬁnition 2A.8]{KST-monoidal-crypto}, is given as follows:
    \begin{equation*}
        \operatorname{gap}_{\K}(\monoid) = \min\{\dim_{\K}(L) \hspace{0.1cm}|\hspace{0.1cm} L \not\cong \onebt \hspace{0.1cm }\text{is a simple $\monoid$ representation} \}.
    \end{equation*}
\end{Theorem}

\begin{proof}
    Given that we have shown that the Motzkin monoid is well connected and has $H^1(\group, \K) \cong 0$, cf. \autoref{P:Split}, since $\monoid$ is well-connected and we have shown that $H^1(\mathcal{M}o_n,\K) \cong 0$, this follows directly from \autoref{RepGap!}.
\end{proof}

\autoref{simple=gap} means that there are no nontrivial extension containing only factors $\onebt$.

\section{Representation Gap of the Motzkin Monoid}\label{S:Main2}

We now come to the main results.

\subsection{Semisimple Representation Gap}
\label{semisimple}
In the following sections, we will prove a lower bound for the semi\-simple representation gap, simple representation gap and the faithful representation gap.

We begin with the semisimple representation gap. This is the most straightforward, defined in \autoref{ssdim} as the dimension of the smallest semisimple non-trivial representation, which is just the sizes of the L-cells. 

Now, given what was observed in \autoref{strandsdiagram}, it makes sense to truncate the monoid to include only large representations. As justified in \autoref{subquotient}, we choose $k=\sqrt{n}$ as our truncated number of through strands. We thus consider the representations for $\mathcal{M}o_n^{\leq \sqrt{n}}$.

\begin{Proposition}
    \label{initial-lowerbounds}
    The dimension of a semisimple representation $L_K$ of $\mathcal{M}o_n^{\leq \sqrt{n}}$ has the following lower bound:
    \[
    \operatorname{ssdim}(L_K) \geq \sum_{t=0}^{n}\frac{\sqrt{n}+1}{\sqrt{n}+t+1} \binom{n}{t} \binom{n-t}{\sqrt{n} + t}.
    \]
\end{Proposition}
\begin{proof}
    From \autoref{motzkin formula}, we have: 
    \[
        |\lcell| = |\rcell| = \sum_{t=0}^{n}\frac{k+1}{k+t+1}\binom{n}{k+2t}\binom{k+2t}{t}.
    \]
    Hence, by \autoref{ssdim}, $\operatorname{ssdim}(L_K) = \sum_{t=0}^{n}\frac{k+1}{k+t+1}\binom{n}{k+2t}\binom{k+2t}{t}.$ Now, consider the representations for $\mathcal{M}o_n^{\leq \sqrt{n}}$. As such, our formula can be rewritten as follows:\begin{equation*}
    \sum_{t=0}^{n}\frac{\sqrt{n}+1}{\sqrt{n}+t+1}\binom{n}{\sqrt{n}+2t}\binom{\sqrt{n}+2t}{t}.
\end{equation*}
This can be expressed more clearly via the following manipulations:
\begin{gather*}
    \sum_{t=0}^{n}\frac{\sqrt{n}+1}{\sqrt{n}+t+1} (\frac{n!(\sqrt{n} + 2t)!}{(n - \sqrt{n} - 2t)! (\sqrt{n} + 2t)!(\sqrt{n} + 2t - t)! t!})  \\
    = \sum_{t=0}^{n}\frac{\sqrt{n}+1}{\sqrt{n}+t+1} (\frac{n!}{(n - t -(\sqrt{n} + t))! (\sqrt{n} + t)! t!}) \\
    = \sum_{t=0}^{n}\frac{\sqrt{n}+1}{\sqrt{n}+t+1} (\frac{n!}{\frac{(n-t)!(n - t -(\sqrt{n} + t))!}{(n-t)!} (\sqrt{n} + t)! t!}) \\
    = \sum_{t=0}^{n}\frac{\sqrt{n}+1}{\sqrt{n}+t+1} (\frac{n!(n-t)!}{(n-t)!(n - t -(\sqrt{n} + t))! (\sqrt{n} + t)! t!}) \\
    = \sum_{t=0}^{n}\frac{\sqrt{n}+1}{\sqrt{n}+t+1} (\frac{n!}{(n-t)!t!}) (\frac{(n-t)!}{(n - t -(\sqrt{n} + t))! (\sqrt{n} + t)!}) \\
    = \sum_{t=0}^{n}\frac{\sqrt{n}+1}{\sqrt{n}+t+1} \binom{n}{t} \binom{n-t}{\sqrt{n} + t}.
\end{gather*}
As proposed.
\end{proof}

We now prove an asymptotic for the growth of $\operatorname{ssgap}(\mathcal{M}o_n^{\leq \sqrt{n}})$ of as $n$ becomes increasingly large. To this end, recall Bachman--Landau notation (``capital O notation''), in particular, $\Theta$ meaning an asymptotic lower and upper bound

\begin{Proposition}
    The semisimple representation gap of $\mathcal{M}o_n^{\leq \sqrt{n}}$ has the following asymptotic lower bound:
    \[
    \operatorname{ssgap}(\mathcal{M}o_n^{\leq \sqrt{n}}) \geq \Theta (n^{-\frac{3}{2}}\cdot 3^n)
    \]
\end{Proposition}
\begin{proof}
    We know from \autoref{initial-lowerbounds} that $\sum_{t=0}^{n}\frac{\sqrt{n}+1}{\sqrt{n}+t+1} \binom{n}{t} \binom{n-t}{\sqrt{n} + t}$ is a lower bound for the size of semisimple representations of the Motzkin mon\-oid over $n$ nodes. Hence, by \autoref{ssdim}, $$\operatorname{ssgap}(\mathcal{M}o_n) \geq \sum_{t=0}^{n}\frac{\sqrt{n}+1}{\sqrt{n}+t+1}  \binom{n}{t} \binom{n-t}{\sqrt{n} + t}.$$ We will now find an asymptotic for this formula.

    \vspace{0.2cm}
    To do so, we wish to rewrite our formula in terms of only $n$, so we find an optimal $t$ with respect to $n$ which we can substitute in. To find this, we plotted the values of our equation over $t \in [0,n]$, and once again found dominant growth coming from a small number of terms, as is clear by  \autoref{t600} below. Note that we change the range of the x-axis to go to 500 since the expression approached 0 as $t$ approached 600:

\begin{figure}[hbt!]
    \centering
    \includegraphics[width=0.45\linewidth]{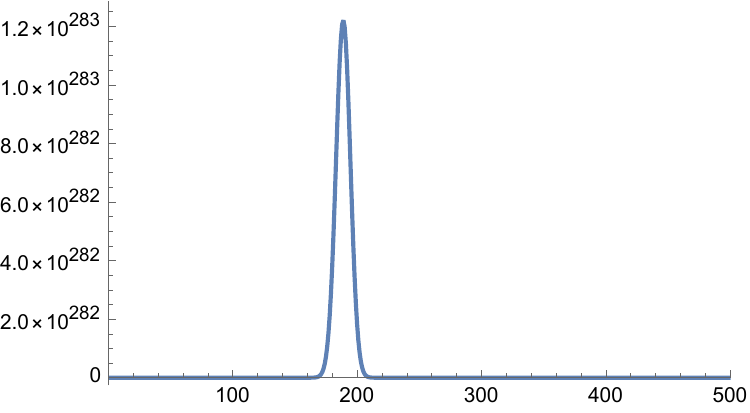}
    \includegraphics[width=0.45\linewidth]{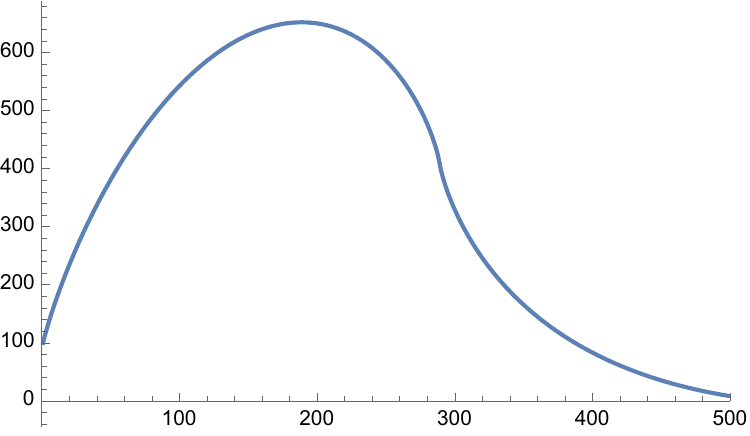}
    \caption{Value within summation of $\operatorname{ssdim}(\mathcal{M}o_{600})$ for varying values of $t \in [0,500]$. The second plot uses a logarithmic scale on the y-axis.}
    \label{t600}
\end{figure}

Note here that the $x$-axis represents values of $t$, and the $y$-axis outputs the values of $\frac{\sqrt{n}+1}{\sqrt{n}+t+1} \binom{n}{t} \binom{n-t}{\sqrt{n} + t}$. In this case, for $n = 600$, there is a clear peak which occurs roughly for $t=200$. We hence assume $t = \frac{n}{3}$ as an approximation for the optimal $t$. (This is not the precise optimal $t$, but good enough to get the lower bound we want.) \autoref{t=n/3} graphs the values of $\frac{n}{3}$, compared with values of $t$ for which $\binom{n}{t}\binom{n-t}{\sqrt{n} + t}$ was the largest.

\begin{figure}[hbt!]
    \centering
    \includegraphics[width=0.45\linewidth]{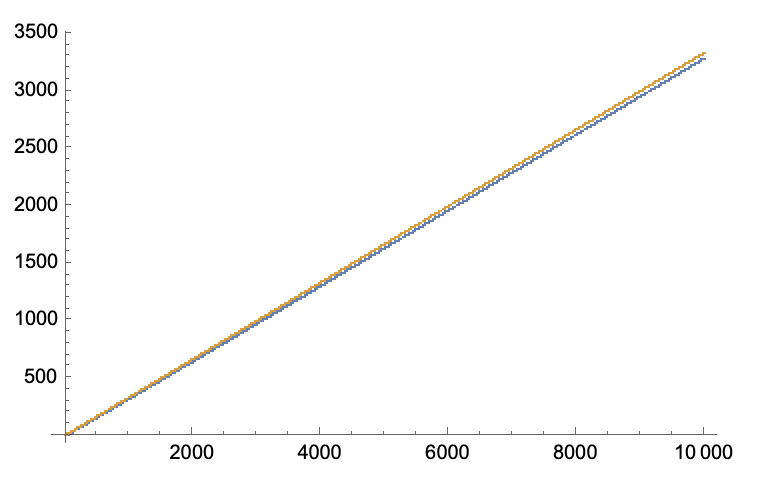}
    \includegraphics[width=0.45\linewidth]{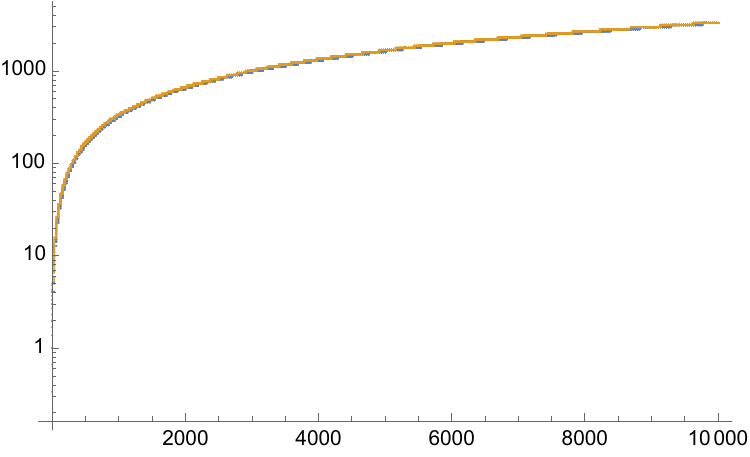}
    \caption{Value of $t$ which yielded the highest output, and the line $t=\frac{n}{3}$. The second plot uses a logarithmic scale on the y-axis.}
    \label{t=n/3}
\end{figure}

We observe that $t = \frac{n}{3}$ consistently gives a close estimation for the maximum term in the sum, which motivated our choice. This allows us to reformulate our calculations with $n$ as the only variable, giving the lower bound; $$\frac{\sqrt{n}+1}{\sqrt{n}+\frac{n}{3}+1} \binom{n}{\frac{n}{3}} \binom{n-\frac{n}{3}}{\sqrt{n} + t} \leq \sum_{t=0}^{n}\frac{\sqrt{n}+1}{\sqrt{n}+t+1} \binom{n}{t} \binom{n-t}{\sqrt{n} + \frac{n}{3}}.$$ From Mathematica, we then find the asymptotic of this to be as follows:
\[
\frac{\sqrt{n}+1}{\sqrt{n}+\frac{n}{3}+1} \binom{n}{\frac{n}{3}} \binom{n-\frac{n}{3}}{\sqrt{n} + \frac{n}{3}} \sim_{n\rightarrow \infty} n^{-\frac{3}{2}} \cdot 3^n, 
\]
so that we obtain,
\[
\operatorname{ssgap}(\mathcal{M}o_n^{\leq \sqrt{n}}) \geq \Theta (n^{-\frac{3}{2}}\cdot 3^n),
\]
as a lower bound for the semisimple representation gap.
\end{proof}

\subsection{Representation Gap}
\label{simple}
We consider here the dimension of representations using Gram matrices as defined earlier. We have that $\mathcal{H}(e)$ is the trivial group, and thus elements in $\K \mathcal{H}(e)$ can be identified as elements in $\K$. Thus, we construct $P^\mathcal{J}$  with reference to any given Green's cell $\jcell$ as follows:
\[
[P^\mathcal{J}_{mn}] = \begin{cases}
    1 \hspace{0.1cm} \text{if } \hspace{0.1cm} [\hcell_{mn}] \hspace{0.1cm} \text{is idempotent}, \\
    0\hspace{0.1cm} \text{otherwise.}
\end{cases}
\]
The rank of this matrix is the dimension of the associated simple representation. Note that this rank depends on the ground field $\K$.

\begin{Example}
    In order to exemplify the construction of a relevant Gram matrix, recall the Motzkin diagrams for $n=3$, specifically we focus on $\jcell_1 \in \mathcal{M}o_3$, again drawing in red the diagrams which are idempotent:
    \begin{gather*}
    \begin{tabular}{C!{\color{blue}\vrule width 0.1mm}C!{\color{blue}\vrule width 0.1mm}C!{\color{blue}\vrule width 0.1mm}C!{\color{blue}\vrule width 0.1mm}C}
\arrayrulecolor{blue}
\begin{tikzpicture}[anchorbase]
\draw[ red,ultra thick] (0,0.75) to (0,0);
\draw[ red,ultra thick] (0.5,0.75) to[out=270,in=180] (0.75,0.5) to[out=0,in=270] (1,0.75);
\draw[ red,ultra thick] (0.5,0) to[out=90,in=180] (0.75,0.25) to[out=0,in=90] (1,0);
\end{tikzpicture} & 
\begin{tikzpicture}[anchorbase]
\draw[ red,ultra thick] (0,0.75) to (0,0);
\draw[ red,ultra thick] (0.5,0.75) to[out=270,in=180] (0.75,0.5) to[out=0,in=270] (1,0.75);
\draw[ red,ultra thick] (0.5,0) to (0.5,0.25);
\draw[ red,ultra thick] (1,0) to (1,0.25);
\end{tikzpicture} &
\begin{tikzpicture}[anchorbase]
\draw[usual] (0,0.75) to (0.5,0);
\draw[usual] (0.5,0.75) to[out=270,in=180] (0.75,0.5) to[out=0,in=270] (1,0.75);
\draw[usual] (0,0) to (0,0.25);
\draw[usual] (1,0) to (1,0.25);
\end{tikzpicture} &
\begin{tikzpicture}[anchorbase]
\draw[ red,ultra thick] (0,0.75) to (1,0);
\draw[ red,ultra thick] (0.5,0.75) to[out=270,in=180] (0.75,0.5) to[out=0,in=270] (1,0.75);
\draw[ red,ultra thick] (0,0) to (0,0.2);
\draw[ red,ultra thick] (0.5,0) to (0.5,0.2);
\end{tikzpicture} &
\begin{tikzpicture}[anchorbase]
\draw[usual] (0,0.75) to (1,0);
\draw[usual] (0.5,0.75) to[out=270,in=180] (0.75,0.5) to[out=0,in=270] (1,0.75);
\draw[usual] (0,0) to[out=90,in=180] (0.25,0.25) to[out=0,in=90] (0.5,0);
\end{tikzpicture} 
\\
\hline
\begin{tikzpicture}[anchorbase]
\draw[ red,ultra thick] (0,0.75) to (0,0);
\draw[ red,ultra thick] (0.5,0.75) to (0.5,0.5);
\draw[ red,ultra thick] (1,0.75) to (1,0.5);
\draw[ red,ultra thick] (0.5,0) to[out=90,in=180] (0.75,0.25) to[out=0,in=90] (1,0);
\end{tikzpicture} & 
\begin{tikzpicture}[anchorbase]
\draw[ red,ultra thick] (0,0.75) to (0,0);
\draw[ red,ultra thick] (0.5,0.75) to (0.5,0.5);
\draw[ red,ultra thick] (1,0.75) to (1,0.5);
\draw[ red,ultra thick] (0.5,0) to (0.5,0.25);
\draw[ red,ultra thick] (1,0) to (1,0.25);
\end{tikzpicture} &
\begin{tikzpicture}[anchorbase]
\draw[usual] (0,0.75) to (0.5,0);
\draw[usual] (0.5,0.75) to (0.5,0.5);
\draw[usual] (1,0.75) to (1,0.5);
\draw[usual] (0,0) to (0,0.25);
\draw[usual] (1,0) to (1,0.25);
\end{tikzpicture} &
\begin{tikzpicture}[anchorbase]
\draw[usual] (0,0.75) to (1,0);
\draw[usual] (0.5,0.75) to (0.5,0.55);
\draw[usual] (1,0.75) to (1,0.55);
\draw[usual] (0,0) to (0,0.2);
\draw[usual] (0.5,0) to (0.5,0.2);
\end{tikzpicture} &
\begin{tikzpicture}[anchorbase]
\draw[usual] (0,0.75) to (1,0);
\draw[usual] (0.5,0.75) to (0.5,0.55);
\draw[usual] (1,0.75) to (1,0.55);
\draw[usual] (0,0) to[out=90,in=180] (0.25,0.25) to[out=0,in=90] (0.5,0);
\end{tikzpicture} 
\\
\hline
\begin{tikzpicture}[anchorbase]
\draw[usual] (0.5,0.75) to (0,0);
\draw[usual] (0,0.75) to (0,0.5);
\draw[usual] (1,0.75) to (1,0.5);
\draw[usual] (0.5,0) to[out=90,in=180] (0.75,0.25) to[out=0,in=90] (1,0);
\end{tikzpicture} & 
\begin{tikzpicture}[anchorbase]
\draw[usual] (0.5,0.75) to (0,0);
\draw[usual] (0,0.75) to (0,0.5);
\draw[usual] (1,0.75) to (1,0.5);
\draw[usual] (0.5,0) to (0.5,0.25);
\draw[usual] (1,0) to (1,0.25);
\end{tikzpicture} &
\begin{tikzpicture}[anchorbase]
\draw[ red,ultra thick] (0.5,0.75) to (0.5,0);
\draw[ red,ultra thick] (0,0.75) to (0,0.5);
\draw[ red,ultra thick] (1,0.75) to (1,0.5);
\draw[ red,ultra thick] (0,0) to (0,0.25);
\draw[ red,ultra thick] (1,0) to (1,0.25);
\end{tikzpicture} &
\begin{tikzpicture}[anchorbase]
\draw[usual] (0.5,0.75) to (1,0);
\draw[usual] (0,0.75) to (0,0.55);
\draw[usual] (1,0.75) to (1,0.55);
\draw[usual] (0,0) to (0,0.2);
\draw[usual] (0.5,0) to (0.5,0.2);
\end{tikzpicture} &
\begin{tikzpicture}[anchorbase]
\draw[usual] (0.5,0.75) to (1,0);
\draw[usual] (0,0.75) to (0,0.55);
\draw[usual] (1,0.75) to (1,0.55);
\draw[usual] (0,0) to[out=90,in=180] (0.25,0.25) to[out=0,in=90] (0.5,0);
\end{tikzpicture}
\\
\hline
\begin{tikzpicture}[anchorbase]
\draw[ red,ultra thick] (1,0.75) to (0,0);
\draw[ red,ultra thick] (0,0.75) to (0,0.55);
\draw[ red,ultra thick] (0.5,0.75) to (0.5,0.55);
\draw[ red,ultra thick] (0.5,0) to[out=90,in=180] (0.75,0.25) to[out=0,in=90] (1,0);
\end{tikzpicture} & 
\begin{tikzpicture}[anchorbase]
\draw[usual] (1,0.75) to (0,0);
\draw[usual] (0.5,0.75) to (0.5,0.55);
\draw[usual] (0,0.75) to (0,0.55);
\draw[usual] (0.5,0) to (0.5,0.2);
\draw[usual] (1,0) to (1,0.2);
\end{tikzpicture} &
\begin{tikzpicture}[anchorbase]
\draw[usual] (1,0.75) to (0.5,0);
\draw[usual] (0,0.75) to (0,0.5);
\draw[usual] (0.5,0.75) to (0.5,0.5);
\draw[usual] (0,0) to (0,0.25);
\draw[usual] (1,0) to (1,0.25);
\end{tikzpicture} &
\begin{tikzpicture}[anchorbase]
\draw[ red,ultra thick] (1,0.75) to (1,0);
\draw[ red,ultra thick] (0,0.75) to (0,0.55);
\draw[ red,ultra thick] (0.5,0.75) to (0.5,0.55);
\draw[ red,ultra thick] (0,0) to (0,0.2);
\draw[ red,ultra thick] (0.5,0) to (0.5,0.2);
\end{tikzpicture} &
\begin{tikzpicture}[anchorbase]
\draw[ red,ultra thick] (1,0.75) to (1,0);
\draw[ red,ultra thick] (0,0.75) to (0,0.55);
\draw[ red,ultra thick] (0.5,0.75) to (0.5,0.55);
\draw[ red,ultra thick] (0,0) to[out=90,in=180] (0.25,0.25) to[out=0,in=90] (0.5,0);
\end{tikzpicture}
\\
\hline
\begin{tikzpicture}[anchorbase]
\draw[usual] (1,0.75) to (0,0);
\draw[usual] (0,0.75) to[out=270,in=180] (0.25,0.5) to[out=0,in=270] (0.5,0.75);
\draw[usual] (0.5,0) to[out=90,in=180] (0.75,0.25) to[out=0,in=90] (1,0);
\end{tikzpicture} & 
\begin{tikzpicture}[anchorbase]
\draw[usual] (1,0.75) to (0,0);
\draw[usual] (0,0.75) to[out=270,in=180] (0.25,0.5) to[out=0,in=270] (0.5,0.75);
\draw[usual] (0.5,0) to (0.5,0.25);
\draw[usual] (1,0) to (1,0.25);
\end{tikzpicture} &
\begin{tikzpicture}[anchorbase]
\draw[usual] (1,0.75) to (0.5,0);
\draw[usual] (0,0.75) to[out=270,in=180] (0.25,0.5) to[out=0,in=270] (0.5,0.75);
\draw[usual] (0,0) to (0,0.25);
\draw[usual] (1,0) to (1,0.25);
\end{tikzpicture} &
\begin{tikzpicture}[anchorbase]
\draw[ red,ultra thick] (1,0.75) to (1,0);
\draw[ red,ultra thick] (0,0.75) to[out=270,in=180] (0.25,0.5) to[out=0,in=270] (0.5,0.75);
\draw[ red,ultra thick] (0,0) to (0,0.2);
\draw[ red,ultra thick] (0.5,0) to (0.5,0.2);
\end{tikzpicture} &
\begin{tikzpicture}[anchorbase]
\draw[ red,ultra thick] (1,0.75) to (1,0);
\draw[ red,ultra thick] (0,0.75) to[out=270,in=180] (0.25,0.5) to[out=0,in=270] (0.5,0.75);
\draw[ red,ultra thick] (0,0) to[out=90,in=180] (0.25,0.25) to[out=0,in=90] (0.5,0);
\end{tikzpicture}
\end{tabular}
    \end{gather*}
    Now, we can draw a matrix where red (idempotent) diagrams correspond to entries of one, and black (non-idempotent) diagrams correspond to entries of zero;
    \begin{center}
    \[\begin{bmatrix}
    1 & 1 & 0&1&0\\
    1 & 1 & 0&0&0\\
    0 & 0 & 1&0&0\\
    1 & 0 & 0&1&1\\
    0 & 0 & 0&1&1\\
    \end{bmatrix}
    \]
    \end{center}
    and this is exactly $P^\mathcal{J}$ for $\mathcal{M}o_3$ for $k=1$ through strands. This matrix is full rank (over any field), so the dimension of the associated simple representation is five (over any field).
\end{Example}

From \autoref{simple=gap}, the simple representations are the smallest possible non-trivial representations of our monoid, so that the simple representation gap is also the general representation gap of our monoid. Referring to \autoref{rank}, we know that the simple dimensions $L_K$ can be calculated as follows $\dim_\K(L_K) = \operatorname{rank}(P^\mathcal{J})$. Hence, we analyse the rank of the Gram matrices associated with our $\jcell$'s in order to obtain a lower bound for the simple representation gap. To do so, we will find full rank sub matrices which determine a minimum value for the rank, hence obtaining a lower bound for the growth of the simple dimensions.

We now prove the claim of this section. In order to analyse how the rank of this matrix grows as we increase $n$, we consider the structure of idempotent diagrams in our monoid.

Recall that an idempotent diagram $a$ is one which remains identical when composed with itself, i.e., $s \cdot s = s$. We note that in order for this to be true, each through strand in the composition must begin and end at the same nodes as it does in the original $a$. This may seem obvious, but it is an important fact in proving the following proposition.

\begin{Proposition}
    \label{idempotency}
    Idempotent diagrams $s$ occur when diagrams have the following local configurations of through strands:
    \begin{enumerate}
        \item Straight through strands, which connect top node $i$ with bottom node $i$.
        \item Curved through strands with cap-cup chains; if the strand connects to node $a$ on top and node $b$ on the bottom, for  $ a < c_i< b$, there exists a sequence of cups and caps;
        \[
        (a,c_1),(c_1,c_2),(c_2,c_3),...,(c_n,b)
        \]
        where $(a,c_1)$ is a cap,  $(c_n,b)$ is a cup, and $(c_i,c_{i+1})$ is a cup for even $i$, and a cap for odd $i$. The following is an example.

    \begin{gather*}
    \begin{tikzpicture}[anchorbase]
    \draw[usual] (6.7,1.5)to[out=0,in=90](7,1);
    \draw[usual] (6.7,1.5) to (6.8,1.3);
    \draw[usual] (6.7,1.5) to (6.9,1.65);
    \draw[usual] (0,1) to[out=270,in=180] (2.5,0.3)to[out=360,in=160](7,0);
    \draw[thick, blue] (7,0) node[right]{$b$};
    \draw[thick, blue] (7,1) node[right]{$b$};
    \draw[thick, blue] (0,1) node[left]{$a$};
    \draw[usual] (0.5,1) to[out=270,in=180] (0.75,0.5) to[out=0,in=270] (1,1);
    \draw[thick, blue] (0.6,1) to (0.8,1.7)node[right]{$c_1$};
    \draw[usual] (0,1) to[out=90,in=180] (0.25,1.5) to[out=0,in=90] (0.5,1);
    \draw[usual] (1,1) to[out=90,in=180] (1.25,1.5);
    \draw[usual] (2.75,1.5) to[out=0,in=90] (3,1);
    \draw[thick, blue](3.1,1) to (3.3,1.7)node[above,right]{$c_{i-2}$};
    \draw[usual] (2.5,1) to[out=270,in=180] (2.75,0.75) to[out=0,in=270] (3,1);
    \draw[thick, blue](2.6,1) to (2.6,1.7)node[above]{$c_{i-1}$};
    \draw[usual] (2,1) to[out=90,in=180] (2.25,1.25) to[out=0,in=90] (2.5,1);
    \draw[thick, blue](1.9,1) to (1.7,0)node[left]{$c_{i}$};
    \draw[usual] (2,1) to[out=270,in=180] (2.75,0.5) to[out=0,in=270] (3.5,1);
    \draw[thick, blue](3.6,1) to (3.8,0)node[right]{$c_{i+1}$};
    \draw[usual] (2,1.5) to (2.75,1.5);
    \draw[usual] (1.5,1.5) to (1.53,1.5);
    \draw[usual] (1.7,1.5) to (1.73,1.5);
    \draw[usual] (1.6,1.5) to (1.63,1.5);
    \draw[usual] (6.5,1) to[out=270,in=180] (6.75,0.5) to[out=0,in=270] (7,1);
    \draw[thick, blue] (6.4,1) to (6,0)node[left]{$c_n$};
    \draw[usual] (6,1) to[out=90,in=180] (6.25,1.5) to[out=0,in=90] (6.5,1);
    \draw[usual] (4.35,1) to (4.4,1);
    \draw[usual] (4.7,1) to (4.75,1);
    \draw[usual] (5.05,1) to (5.1,1);
    \draw[usual] (3.5,1.1) to (3.5,1.07);
    \draw[usual] (3.55,1.2) to (3.55,1.17);
    \draw[usual] (3.65,1.3) to (3.65,1.27);
    \draw[usual] (5.85,0.7) to (5.85,0.73);
    \draw[usual] (5.95,0.8) to (5.95,0.83);
    \draw[usual] (6,0.9) to (6,0.93);
    \end{tikzpicture}
    \end{gather*}
    \end{enumerate}

\end{Proposition}

\begin{proof}
    In order for $s$ to be idempotent, it must retain the same number of through strands. For this to be possible, the top of each through strand must be linked with its bottom in the composition. 
    It is clear that (1) is true since straight through strands will match up with themselves by construction and thus are always idempotent. 
    If a through strand is not straight, it must somehow be linked to itself via composition. This happens via the construction in (2), whereby a through strand connecting $a$ to $b$ is preserved under composition, with $s \cdot s$ having a through strand from $a$ to $b$, ensuring that $s$ is locally idempotent. We will prove, by contradiction, that this construction is necessary. Suppose there is an idempotent cell which does not have caps and cups in the above defined sequence. This could happen in four cases, for each of which we give a reason why this would result in loss of idempotency:
    \begin{enumerate}
        \item There are no caps connected to $a$:
    
            If there is not a cap connecting to $a$, it must be a disconnected node or connect to a through strand. If it is a disconnected node, then the through strand will attached to $b$ becomes disconnected in composition and there is a clear contradiction. If it has a through strand, then there exists in our diagram some strand which goes from some $c$ to $a$, so in $s \cdot s$, we will have a through strand from $c$ to $b$, a strand which was not in the original $s$, meaning that we are left with a different diagram under composition, and $s\cdot s \not= s$.
            \item There are no cups connected to $b$:
    
            The proof of this is equivalent to the proof of above, replacing caps with cups.
    
            \item There is a gap in the sequence such that $(c_i,c_{i+1})$ is followed by a cap or cup beginning at $c_{i+n}$ for $n \geq 2$:

            This will break our construction, as if we have a sequence of connected cups and caps, and then it stops at $c_i$ where there is either a through strand or disconnected node, we have a situation equivalent to the first case, where the through strand becomes disconnected or changes its endpoints, violating idempotency.
    
            \item         $(c_i,c_{i+1})$ is a cap for even $i$ or a cup for odd $i$:
            
            This also breaks our construction. If the caps and cups are consecutive, this would imply that there is are two caps or cups coming from the same node, which is impossible by construction of the Motzkin monoid.
    \end{enumerate}
    
    Hence, if a diagram with a non-straight through strand fails to meet the requirements of (2), it is not idempotent, proving that (2) is the only suitable local construction of an idempotent diagram with non-straight through strands.
\end{proof}

Now, since our rank is determined by the number of linearly independent rows or columns, we wish to consider rows that will always be linearly independent. One can observe that symmetry in through strands can cause row reductions (see \autoref{rowreductions} below for an instance of this). We hence consider rows which consist of only diagrams with all through strands pushed to one side at the top, and proving that these will be linearly independent. We begin by showing the following:

\begin{Proposition}
    \label{pairwise lin ind}
    Any two distinct $\rcell$'s or rows with top presentations having all through strands pushed to one side will be linearly independent.
\end{Proposition}

\begin{proof}
    We initially make clear what this presentation looks like. For this proof, we assume strands are all on the left; the right side argument is equivalent. Consider a diagram with all strands pushed to the left, with the final strand starting at node $a$ and all nodes to the right of $a$ either connected to a cup or disconnected. For example, if we let $a$ be the third node we could have the following as a possible top row presentation: 

    \begin{gather*}
    \begin{tikzpicture}[anchorbase]
    \draw[usual] (0,0) to (0,1);
    \draw[usual] (0.5,0) to (0.5,1);
    \draw[usual] (1,1) to (1,0);
    \draw[usual] (1.5,1) to[out=270,in=180] (1.75,0.75) to[out=0,in=270] (2,1);
    \draw[usual] (2.5,1) to (2.5,0.75);
    \draw[usual] (3,1) to[out=270,in=180] (4,0.5) to[out=0,in=270] (5,1);
    \draw[usual] (3.5,1) to (3.5,0.75);
    \draw[usual] (4,1) to[out=270,in=180] (4.25,0.75) to[out=0,in=270] (4.5,1);
    \draw[usual] (5.5,1) to[out=270,in=180] (5.75,0.75) to[out=0,in=270] (6,1);
    \end{tikzpicture}
    \end{gather*}
    
    From above, the cups on the top exist between pairs of nodes $(i_k,j_{k'})$, for $a < i_k, j_{k'} \leq n$. Suppose we have a top row with $k$ cups, we can define its cup set: $T := \{ (i_k,j_{k'}), (i_l,b) \}$, where $i_k$ indicates the node which each cup begins, and $j_{k'}$, the node where the cup ends. We use $b$ to denote the node at which the final cup terminates.
    Then, by our above definition, the bottom presentations which form idempotent diagrams in combination with this top row will have caps at, 
    \begin{gather*}B = \{(a,\lambda_{a}),(\lambda_b,\lambda_{c})| \text{where } \lambda_{a,b,c,\cdots} \in T  \\
    \text{represent arbitrary choices of } i_k \text{ or } j_{k'} \text{ with } (i_k,j_{k'}) \in T, \\ \text{ filling the requirements as in \autoref{idempotency}(b)}, \}\end{gather*} ensuring that the through strand $(a,b)$ is locally idempotent. Note that in the definition of the bottom cap set, multiple different combinations of $(\lambda_b,\lambda_{c})$ can be chosen such that we will have a set of possible sets, and there is not one unique idempotent presentation. We assume the bottom presentation to have a through strand at $b$, the end of the final cup, which will be connected to $a$. This means that the bottom presentation of the idempotent elements is wholly determined by the presentation of the top row. Before continuation of the proof, we give an example to clarify what is meant by top and bottom sets.
    
    \begin{Example}
        The following:
        \begin{gather*}
        \begin{tikzpicture}[anchorbase]
    \draw[usual] (0,0.5) node[above]{$a$} to (0,0);
    \draw[usual] (0.5,0.5) node[above]{$1$} to[out=270,in=180] (1,0.25) to[out=0,in=270] (1.5,0.5) node[above]{$3$};
    \draw[usual] (1,0.5)node[above]{$2$} to (1,0.35);
    \draw[usual] (2,0.5) node[above]{$4$} to[out=270,in=180] (2.75,0.2) to[out=0,in=270] (3.5,0.5) node[above]{$7$};
    \draw[usual] (2.5,0.5) node[above]{$5$} to[out=270,in=180] (2.75,0.35) to[out=0,in=270] (3,0.5) node[above]{$6$};
    \draw[usual] (4,0.5) node[above]{$8$} to[out=270,in=180] (4.25,0.25) to[out=0,in=270] (4.5,0.5) node[above]{$b$};
    \end{tikzpicture}
        \end{gather*}
        would have cup set: 
        \[
        \{ (1,3),(5,6),(4,7),(8,b) \}.
        \]
        Consequently the cap set would include, but not be limited to, the following:
        \[
        \{(a,1),(3,6),(4,5),(7,8)\},\{(a,1),(3,8)\}, \{(a,8)\},...
        \]    
        as is diagrammatically clear from the following images:
        \begin{gather*}
        \begin{tikzpicture}[anchorbase]
    \draw[usual] (0,0.5) to (0,0);
    \draw[red,ultra thick] (0,0.5) to[out=90,in=180] (0.25,0.75) to[out=0,in=90] (0.5,0.5);
    \draw[usual] (0.5,0.5) to[out=270,in=180] (1,0.25) to[out=0,in=270] (1.5,0.5);
    \draw[red,ultra thick] (1.5,0.5) to[out=90,in=180] (2.25,0.8) to[out=0,in=90] (3,0.5);
    \draw[usual] (1,0.5) to (1,0.35);
    \draw[usual] (2,0.5) to[out=270,in=180] (2.75,0.2) to[out=0,in=270] (3.5,0.5);
    \draw[red,ultra thick] (2,0.5) to[out=90,in=180] (2.25,0.65) to[out=0,in=90] (2.5,0.5);
    \draw[usual] (2.5,0.5)  to[out=270,in=180] (2.75,0.35) to[out=0,in=270] (3,0.5);
    \draw[usual] (4,0.5)  to[out=270,in=180] (4.25,0.25) to[out=0,in=270] (4.5,0.5);
    \draw[red,ultra thick] (3.5,0.5) to[out=90,in=180] (3.75,0.65) to[out=0,in=90] (4,0.5);
    \draw[usual] (4.5,0.5) to (4.5,1);
    \end{tikzpicture}
    \hspace{1cm}
    \begin{tikzpicture}[anchorbase]
    \draw[usual] (0,0.5) to (0,0);
    \draw[red,ultra thick] (0,0.5) to[out=90,in=180] (0.25,1) to[out=0,in=90] (0.5,0.5);
    \draw[usual] (0.5,0.5) to[out=270,in=180] (1,0.25) to[out=0,in=270] (1.5,0.5);
    \draw[usual] (1,0.5) to (1,0.35);
    \draw[usual] (2,0.5) to[out=270,in=180] (2.75,0.2) to[out=0,in=270] (3.5,0.5);
    \draw[usual] (2.5,0.5)  to[out=270,in=180] (2.75,0.35) to[out=0,in=270] (3,0.5);
    \draw[red,ultra thick] (1.5,0.5) to[out=90,in=180] (2.75,1) to[out=0,in=90] (4,0.5);
    \draw[usual] (4,0.5)  to[out=270,in=180] (4.25,0.25) to[out=0,in=270] (4.5,0.5);
    \draw[usual] (4.5,0.5) to (4.5,1);
    \end{tikzpicture}
    \hspace{1cm}
    \begin{tikzpicture}[anchorbase]
    \draw[usual] (0,0.5) to (0,0);
    \draw[red,ultra thick] (0,0.5) to[out=90,in=180] (2,1.25) to[out=0,in=90] (4,0.5);
    \draw[usual] (0.5,0.5) to[out=270,in=180] (1,0.25) to[out=0,in=270] (1.5,0.5);
    \draw[usual] (1,0.5) to (1,0.35);
    \draw[usual] (2,0.5) to[out=270,in=180] (2.75,0.2) to[out=0,in=270] (3.5,0.5);
    \draw[usual] (2.5,0.5)  to[out=270,in=180] (2.75,0.35) to[out=0,in=270] (3,0.5);
    \draw[usual] (4,0.5)  to[out=270,in=180] (4.25,0.25) to[out=0,in=270] (4.5,0.5);
    \draw[usual] (4.5,0.5) to (4.5,1);
    \end{tikzpicture}
        \end{gather*}
        This is a rather prototypical example.
    \end{Example}  
    We now prove that each top set uniquely determines a bottom set. Consider two cells $A, A'$ which have cups in the same places; 
    \begin{gather*}
    \{ (i_1,j_{1'}), ..., (i_k,j_{k'})\} = \{ (i_1',j_{1'}'), ..., (i_k',j_{k'}')\}.
    \end{gather*}
    We can thus construct the bottom idempotent sets as defined above to be $B$ and $B'$. However since we assumed the top sets to be the same, we know that $i_l=i_l'$ and $j_{l'}=j_{l'}'$, $\forall l \in [1,k]$. Since $B$ and $B'$ are constructed wholly from these $(i_k,j_{k'})$, we have $B=B'$ which implies that the bottom presentations have equivalent caps, showing that the set of top row cups uniquely determines the set of caps in the bottom presentations. Note that if any cup in $A$ or $A'$ was placed differently, its respective bottom set would be different since the caps would have to be adjusted accordingly to ensure \autoref{idempotency} is satisfied. Thus, unique top presentations have idempotents in combination with different bottom presentations, implying pairwise linear independence.
    
    In this argument, we defined the bottom set to have a strand at $b$, the final cup, because we wanted to prove uniqueness over differing cup sets which all terminate at $b$. Pairwise independence between bottom presentations relating to top presentations which end at a different node, for example node $d$, is obvious, since the bottom presentations will have through strands at $d \neq b$. As such, we ensure that the set of bottom presentations which form idempotent elements are unique for each different top presentation. Now since $\lcell$'s are determined by bottom presentations, idempotent elements in different rows therefore appear in different columns, ensuring that any two rows with unique top presentations will be linearly independent in our Gram matrix. This completes the proof of \autoref{pairwise lin ind}.
\end{proof}

We must now prove linear independence over all rows, not just any combination of two rows. Namely, we prove that no row can be expressed as the linear combination of any other set of rows.

\begin{Example}
    We exemplify this within $\mathcal{M}o_5$ for $\jcell_1$ i.e. the J-cell with one through strand. We fix the top of the through strand at the left most node, and show that it is possible to construct a full-rank sub-cell with respect to the placing of idempotents. We include this sub-cell below, denoting idempotent elements in red.

    \begin{gather*}
\begin{tabular}{C!{\color{blue}\vrule width 0.1mm}C!{\color{blue}\vrule width 0.1mm}C!{\color{blue}\vrule width 0.1mm}C!{\color{blue}\vrule width 0.1mm}C!{\color{blue}\vrule width 0.1mm}C!{\color{blue}\vrule width 0.1mm}C!{\color{blue}\vrule width 0.1mm}C!{\color{blue}\vrule width 0.1mm}C}
\arrayrulecolor{blue}
\begin{tikzpicture}
[anchorbase]
\draw[red, thick] (0,1) to (0,0);
\draw[red, thick] (0.25,1) to[out=270,in=180] (0.375,0.75) to[out=0,in=270] (0.5,1);
\draw[red, thick] (0.75,1) to[out=270,in=180] (0.875,0.75) to[out=0,in=270] (1,1);
\draw[red, thick] (0.25,0) to (0.25,0.25);
\draw[red, thick] (0.5,0) to (0.5,0.25);
\draw[red, thick] (0.75,0) to (0.75,0.25);
\draw[red, thick] (1,0) to (1,0.25);
\end{tikzpicture} & 
\begin{tikzpicture}[anchorbase]
\draw[red, thick] (0,1) to (0.5,0);
\draw[red, thick] (0.25,1) to[out=270,in=180] (0.375,0.75) to[out=0,in=270] (0.5,1);
\draw[red, thick] (0.75,1) to[out=270,in=180] (0.875,0.75) to[out=0,in=270] (1,1);
\draw[red,thick] (0,0) to[out=90,in=180] (0.125,0.25) to[out=0,in=90] (0.25,0);
\draw[red,thick] (0.75,0) to[out=90,in=180] (0.875,0.25) to[out=0,in=90] (1,0);
\end{tikzpicture} &
\begin{tikzpicture}[anchorbase]
\draw[thick] (0,1) to (0.75,0);
\draw[thick] (0.25,1) to[out=270,in=180] (0.375,0.75) to[out=0,in=270] (0.5,1);
\draw[thick] (0.75,1) to[out=270,in=180] (0.875,0.75) to[out=0,in=270] (1,1);
\draw[thick] (0,0) to[out=90,in=180] (0.125,0.25) to[out=0,in=90] (0.25,0);
\draw[thick] (0.5,0) to (0.5,0.2);
\draw[thick] (1,0) to (1,0.25);
\end{tikzpicture} &
\begin{tikzpicture}[anchorbase]
\draw[thick] (0,1) to (0.75,0);
\draw[thick] (0.25,1) to[out=270,in=180] (0.375,0.75) to[out=0,in=270] (0.5,1);
\draw[thick] (0.75,1) to[out=270,in=180] (0.875,0.75) to[out=0,in=270] (1,1);
\draw[thick] (0,0) to[out=90,in=180] (0.25,0.25) to[out=0,in=90] (0.5,0);
\draw[thick] (0.25,0) to (0.25,0.15);
\draw[thick] (1,0) to (1,0.25);
\end{tikzpicture} &
\begin{tikzpicture}[anchorbase]
\draw[thick] (0,1) to (1,0);
\draw[thick] (0.25,1) to[out=270,in=180] (0.375,0.8) to[out=0,in=270] (0.5,1);
\draw[thick] (0.75,1) to[out=270,in=180] (0.875,0.8) to[out=0,in=270] (1,1);
\draw[thick] (0,0) to[out=90,in=180] (0.125,0.25) to[out=0,in=90] (0.25,0);
\draw[thick] (0.5,0) to (0.5,0.2);
\draw[thick] (0.75,0) to (0.75,0.18);
\end{tikzpicture} &
\begin{tikzpicture}[anchorbase]
\draw[thick] (0,1) to (1,0);
\draw[thick] (0.25,1) to[out=270,in=180] (0.375,0.8) to[out=0,in=270] (0.5,1);
\draw[thick] (0.75,1) to[out=270,in=180] (0.875,0.8) to[out=0,in=270] (1,1);
\draw[thick] (0,0) to[out=90,in=180] (0.25,0.25) to[out=0,in=90] (0.5,0);
\draw[thick] (0.25,0) to (0.25,0.18);
\draw[thick] (0.75,0) to (0.75,0.18);
\end{tikzpicture} &
\begin{tikzpicture}[anchorbase]
\draw[red, thick] (0,1) to (1,0);
\draw[red, thick] (0.25,1) to[out=270,in=180] (0.375,0.8) to[out=0,in=270] (0.5,1);
\draw[red, thick] (0.75,1) to[out=270,in=180] (0.875,0.8) to[out=0,in=270] (1,1);
\draw[red,thick] (0.25,0) to[out=90,in=180] (0.375,0.15) to[out=0,in=90] (0.5,0);
\draw[red,thick] (0,0) to[out=90,in=180] (0.375,0.25) to[out=0,in=90] (0.75,0);
\end{tikzpicture} &
\begin{tikzpicture}[anchorbase]
\draw[thick] (0,1) to (1,0);
\draw[thick] (0.25,1) to[out=270,in=180] (0.375,0.8) to[out=0,in=270] (0.5,1);
\draw[thick] (0.75,1) to[out=270,in=180] (0.875,0.8) to[out=0,in=270] (1,1);
\draw[thick] (0,0) to[out=90,in=180] (0.375,0.25) to[out=0,in=90] (0.75,0);
\draw[thick] (0.25,0) to (0.25,0.15);
\draw[thick] (0.5,0) to (0.5,0.15);
\end{tikzpicture} &
\begin{tikzpicture}[anchorbase]
\draw[red, thick] (0,1) to (1,0);
\draw[red, thick] (0.25,1) to[out=270,in=180] (0.375,0.8) to[out=0,in=270] (0.5,1);
\draw[red, thick] (0.75,1) to[out=270,in=180] (0.875,0.8) to[out=0,in=270] (1,1);
\draw[red,thick] (0,0) to[out=90,in=180] (0.125,0.2) to[out=0,in=90] (0.25,0);
\draw[red,thick] (0.5,0) to[out=90,in=180] (0.625,0.2) to[out=0,in=90] (0.75,0);
\end{tikzpicture}
\\
\end{tabular}
    \end{gather*}

    The Gram matrix of this sub-cell can now be written as follows:

\begin{center}
    \[%
    \]
    \end{center}

Clearly, here we have $|\lcell|=9$, and it is not hard to see that the Gram matrix has rank 9, and that we hence have a full rank matrix. Note that this is a sub-matrix of the representation of $\jcell_1 \subseteq \mathcal{M}o_5$, and thus, its rank gives us a lower bound for the entire $\jcell_1$. It is clear that some rows are partial linear combinations of others, sharing elements in the same row, however no column is completely linearly dependent on any other combinations of columns, and the below proof will generalise this result. We also remark that there are other possible bottom representations here, however we do not include them as we consider only linearly independent rows or $\rcell$'s, the number of which is equal to the number of columns.
\end{Example}

\begin{Proposition}
    \label{submonoid}
    There exists a full-rank sub-matrix $M$ of the Gram matrix $P^\jcell$ associated to the diagrammatic $\jcell$ entries over rows which fix all $k$ through strands consecutively on the first $k$ right or left nodes.
\end{Proposition}

\begin{proof}
    What we wish to show here is that for any fixed top row presentation $T_i$, the set of bottom row presentations $B_i$ is unique and not a linear combination of elements from other bottom sets $B_n$ of $T_n \neq T_i$.
    
    We prove this by induction on the number of nodes next to the straight strands in our diagram. To begin, we prove this to be true for three free through strands, since the cases for one and two nodes are trivial. We can show this diagrammatically:
    \begin{gather*}
    \begin{tikzpicture}[anchorbase]
    \draw[blue, ultra thick] (0,0) to (0,0.6);
    \draw[blue, ultra thick] (0.5,0.6) to[out=270,in=180] (0.75,0.3) to[out=0,in=270] (1,0.6);
    \draw[blue, ultra thick] (1.5,0.6) to (1.5,0.35);
    \end{tikzpicture}
    \; \hspace{1cm} \;
    \begin{tikzpicture}[anchorbase]
    \draw[blue, ultra thick] (0,0) to (0,0.6);
    \draw[blue, ultra thick] (1,0.6) to[out=270,in=180] (1.25,0.3) to[out=0,in=270] (1.5,0.6);
    \draw[blue, ultra thick] (0.5,0.6) to (0.5,0.35);
    \end{tikzpicture}
    \; \hspace{1cm} \;
    \begin{tikzpicture}[anchorbase]
    \draw[blue, ultra thick] (0,0) to (0,0.6);
    \draw[blue, ultra thick] (0.5,0.6) to[out=270,in=180] (1,0.3) to[out=0,in=270] (1.5,0.6);
    \draw[blue, ultra thick] (1,0.6) to (1,0.45);
    \end{tikzpicture}
    \end{gather*}
    \begin{gather*}
    \begin{tikzpicture}[anchorbase]
    \draw[ultra thick] (1,0) to (0,1);
    \draw[ultra thick] (0.5,1) to[out=270,in=180] (0.75,0.7) to[out=0,in=270] (1,1);
    \draw[ultra thick] (1.5,1) to (1.5,0.75);
    \draw[ultra thick] (1.5,0) to (1.5,0.25);
    \draw[ultra thick] (0,0) to[out=90,in=180] (0.25,0.25) to[out=0,in=90] (0.5,0);
    \end{tikzpicture}
    \; \hspace{1cm} \;
    \begin{tikzpicture}[anchorbase]
    \draw[ultra thick] (1.5,0) to (0,1);
    \draw[ultra thick] (1,1) to[out=270,in=180] (1.25,0.7) to[out=0,in=270] (1.5,1);
    \draw[ultra thick] (0.5,1) to (0.5,0.85);
    \draw[ultra thick] (0.5,0) to (0.5,0.15);
    \draw[ultra thick] (0,0) to[out=90,in=180] (0.5,0.25) to[out=0,in=90] (1,0);
    \end{tikzpicture}
    \; \hspace{1cm} \;
    \begin{tikzpicture}[anchorbase]
    \draw[ultra thick] (1.5,0) to (0,1);
    \draw[ultra thick] (0.5,1) to[out=270,in=180] (1,0.7) to[out=0,in=270] (1.5,1);
    \draw[ultra thick] (1,1) to (1,0.85);
    \draw[ultra thick] (1,0) to (1,0.15);
    \draw[ultra thick] (0,0) to[out=90,in=180] (0.25,0.25) to[out=0,in=90] (0.5,0);
    \end{tikzpicture}
    \end{gather*}
    \begin{gather*}
    \begin{tikzpicture}[anchorbase]
    \draw[red, ultra thick] (1,0) to (1,0.6);
    \draw[red, ultra thick] (1.5,0) to (1.5,0.25);
    \draw[red, ultra thick] (0,0) to[out=90,in=180] (0.25,0.25) to[out=0,in=90] (0.5,0);
    \end{tikzpicture}
    \; \hspace{1cm} \;
    \begin{tikzpicture}[anchorbase]
    \draw[red, ultra thick] (1.5,0) to (1.5,0.6);
    \draw[red, ultra thick] (0.5,0) to (0.5,0.15);
    \draw[red, ultra thick] (0,0) to[out=90,in=180] (0.5,0.25) to[out=0,in=90] (1,0);
    \end{tikzpicture}
    \; \hspace{1cm} \;
    \begin{tikzpicture}[anchorbase]
    \draw[red, ultra thick] (1.5,0) to (1.5,0.6);
    \draw[red, ultra thick] (1,0) to (1,0.15);
    \draw[red, ultra thick] (0,0) to[out=90,in=180] (0.25,0.25) to[out=0,in=90] (0.5,0);
    \end{tikzpicture}
    \end{gather*}
    Above, the blue row indicates all possible unique top compositions (corresponding to different rows in our $\jcell$), the black represents the idempotent diagrams formed using them, and the red sketches represent the correspondingly unique bottom row idempotent presentations. These will result in entries of one in different columns of the Gram matrix. Since all different rows have entries of one in different columns, we will have a full-rank matrix.

    Now we assume this to be true for $n$ dots and prove it for $n+1$ under this assumption. Here, we only need to consider independence of diagrams which have a cup ending at node $n+1$ since all other diagrams will be linearly independent due to our inductive assumption. Now, the set of diagrams which have a cup ending at $n+1$ will have a subset consisting of idempotent diagrams which do not have a through strand ending at $n+1$. Elements of this subset will be linearly independent as a result of our inductive assumption. 
    
    Hence, we must only prove why the diagrams over $n+1$ nodes will have a set of idempotents which is linearly independent from the set of idempotents on $n$ nodes and within itself. The first is obvious, since none of the idempotents on $n$ strands contain strands going to $n+1$, and thus the latter will consist of elements which do not exist in the former, and hence could not be obtained by linear combinations. 
    
    To prove the second, we consider the sets $T_i, B_i$ and $T_j, B_j$. It may be that $T_i \not\subset T_j$ and $T_j \not\subset T_i$, and it is impossible to construct $T_i + T_j := T_{ij} := \{t|t\in T_i \hspace{0.1cm} \text{or} \hspace{0.1cm} T_j \}$ due to the crossingless property of the Motzkin monoid. For example, the following two top cells could never be equal through the addition of cups:
    \begin{gather*}
    \begin{tikzpicture}[anchorbase]
    \draw[thick] (0.5,1) to (0.5,0.85);
    \draw[thick] (0,1) to[out=270,in=180] (0.5,0.75) to[out=0,in=270] (1,1);
    \draw[thick] (1.5,1) to (1.5,0.85);
    \end{tikzpicture}
    \; \hspace{0.2cm} , \hspace{0.2cm} \;
    \begin{tikzpicture}[anchorbase]
    \draw[thick] (0,1) to (0,0.85);
    \draw[thick] (0.5,1) to[out=270,in=180] (1,0.75) to[out=0,in=270] (1.5,1);
    \draw[thick] (1,1) to (1,0.85);
    \end{tikzpicture}
    \; \hspace{0.2cm} i.e. \hspace{0.2cm} \;
    \begin{tikzpicture}[anchorbase]
    \draw[thick] (0.5,1) to[out=270,in=180] (1,0.75) to[out=0,in=270] (1.5,1);
    \draw[thick] (0,1) to[out=270,in=180] (0.5,0.75) to[out=0,in=270] (1,1);
    \end{tikzpicture}
    \; \hspace{0.2cm} \not\in \mathcal{M}o_4
    \end{gather*}
    In this case, $B_i,B_j$ are linearly independent by \autoref{pairwise lin ind}, and their linear combination is not an element of the monoid, so they do not threaten the possibility for linear independence. As such, we consider only why when $\cup \{ T_1,  \cdots  T_m \} = T_{\sum_m} \in \mathcal{M}o_n$ where $T_i = T_j \Leftrightarrow i = j,$ for $i,j \in [1,m]$, wishing to show that $B_{\sum_m} \neq \cup \{ B_1,  \cdots  B_m \}$. We do this by claiming that there will be a unique idempotent element $\iota \in B_{\sum_m}$ where $\iota \not\in \cup \{ B_1,  \cdots  B_m \}$. This is formed by drawing a bottom presentation with caps connected to every cup in the top presentation, at least one of which was missing before the sum was taken since we are summing over unique elements. $\iota$ has caps with the following construction, the first cap starting at $a$, then, going from left to right across the top:
    \begin{itemize}
        \item If the next cup $(i_k,j_{k'})$ has no cups within it, we end our cap at $i_k$ and begin the next one at $j_{k'}$.
        \item If the next cup has cups within it, we call it a \emph{boundary cup}, and we go to the outermost cup embedded within it $(i_k,j_{k'})$, namely the cup with the left-most $i_k$, ending our cap at $j_{k'}$ to connect to cup $(i_k,j_{k'})$, and begin the next cap at $i_k$. Note in the case of embedded cups, that if there are no more cups to the right, we connect our next cap to $i_k$ in the boundary cup.
    \end{itemize}
    We continue to follow these steps until we end up connecting a cap to the final cup $(i_k,n+1)$, here having created an idempotent diagram. As such, we know that our diagram will have a connected through strand under composition since, as in the construction above, it is generated by the following:
    \begin{itemize}
        \item Consecutive cups, for which we simply join the endpoint of one to the start point of another with a cap.
        \item Embedded cups, for which we begin by connecting to the outer-most embedded cup, and then continue on the next layer, until reaching the innermost cup, which we connect to the beginning of the boundary cup, allowing our through strand to continue out the other end.
    \end{itemize} 
    This construction will always exist as a crossingless presentation since idempotent elements can be built via the above generators using crossingless caps.

    $\iota \not\in \cup \{ B_1,  \cdots  B_m \}$ because of the following argument: Since caps in $\iota$ connect to every cup in $T_{\sum_m}$, the removal of any cup to get some $T_i \in \cup T_m$ will result in the above construction becoming disconnected. Hence, it would no longer be idempotent, thus, $\forall i$, we know that $\iota \not\in B_i$, and thus not in $\cup_i B_i$, showing that $ B_{\sum_m} \neq \cup \{ B_1,  \cdots  B_m \}$ as required. Hence, the inductive hypothesis has been proven true for $n+1$ elements, and thus all $n$. As such, we have proven that it is possible to form a linearly independent set of idempotent bottom presentations on $n$ strands, equal in size to the number of unique top row presentations. Hence, we have proven the existence of a sub-matrix $M$ of $P^\jcell$ which is of full-rank.
\end{proof}

Note that this proof works for clusters of strands on both the left and right, but we do not multiply by 2 since it will not affect exponential growth. Now, since the rank of a matrix is larger than or equal to the rank of any of its sub-matrices, we state the following:

\begin{Proposition}
  \label{simple formula}
  The size of the full-rank submatrix given by the possible diagrammatic presentations where all through strands sit consecutively on the first $n$ nodes, gives a lower bound for the rank of our representation;
  \[
  \operatorname{gap}_\K(\mathcal{M}o_n) \geq |M| = \sum_{t=0}^{n-\sqrt{n}}\frac{1}{t+1}\binom{n-\sqrt{n}}{2t}\binom{2t}{t}. \hspace{0.1cm}
  \] 
\end{Proposition}

\begin{proof}
    We already have the existence of such a submonoid $M$ from \autoref{submonoid}, and it remains to prove its size. To do so we must count the number of ways one can draw caps on the free nodes to the right of the clustered through strands. Conducting this count is possible through use of the formula for the Motzkin monoid, which already counts the number of possible L or R cells; the number of possible top or bottom row presentations on a Motzkin diagram for a given $k$, number of through strands. We are counting the number of different presentations with no through strands, hence taking $k=0$. We are considering possible presentations over $n-\sqrt{n}$ nodes since we need to subtract the nodes where the consecutive through strands will sit, and count the possible combination only on free nodes. We choose $\sqrt{n}$ because it is the optimal number of through strands, as shown in \autoref{subquotient}. Since $\operatorname{gap}_{\K}(\mathcal{M}o_n)$ is defined exactly as the number of simple representations; linearly independent rows/columns in the Gram matrix, this result follows immediately when plugging in these values.
\end{proof}

\begin{figure}[hbt!]
    \centering
    \includegraphics[width=0.5\linewidth]{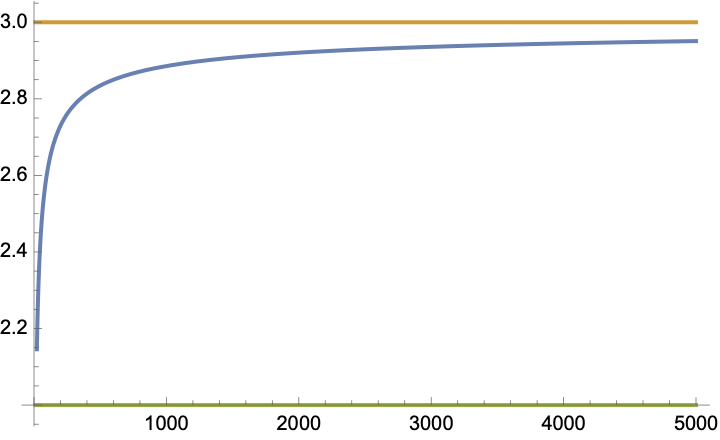}
    \caption{Growth of $n^{th}$ root of our lower bound for $\operatorname{gap}_{\monoid}(\mathcal{M}o_n)$.}
    \label{growthofn}
\end{figure}

We now give a numerical estimate for the growth of this formula, motivated by \autoref{growthofn}, wherein the blue line represents the $n^{th}$ root of our gap, clearly approaching $3$ (but slowly).

We introduce the following lemmas in order to aid in proving \autoref{hypergeo}.
Here and throughout, $2F1(a,b,c;d)$ denotes the ordinary hypergeometric function, as discussed in \cite{Temme-SpecialFunctions}.

\begin{Lemma}
    \label{c1/2}
    The following are equivalent:
    \[
    (a \cdot 2F1(\frac{1}2{(-n+\sqrt{n})},\frac{1}2{(1- n+\sqrt{n})},1;4))^{\frac{1}{n}}
    = (2F1(\frac{1}2{(-n+\sqrt{n})},\frac{1}2{(1- n+\sqrt{n})},\frac{1}{2};4))^{\frac{1}{n}}.
    \]
\end{Lemma}
\begin{proof}
    This can be proven by considering the explicit forms of the hypergeometric functions, these are well known (as in, for example, \cite{Temme-SpecialFunctions}):
    \begin{gather*}
        a \cdot \sum_{n=0}^\infty \frac{(a)_n(b)_n}{(c)_n}\frac{z^n}{n!} = \sum_{n=0}^\infty \frac{(a)_n(b)_n}{(c-\frac{1}{2})_n}\frac{z^n}{n!},
    \end{gather*}
    so for each n;
    \begin{gather*}
        a \cdot \frac{(a)_n(b)_n}{(c)_n}\frac{z^n}{n!} = \frac{(a)_n(b)_n}{(c-\frac{1}{2})_n}\frac{z^n}{n!} \\
        a = \frac{(c)_n}{(c-\frac{1}{2})_n},
    \end{gather*}
    for our case where $c = 1$, $a = \frac{(1)_n}{(\frac{1}{2})_n} \xrightarrow{n\rightarrow \infty} 1$. Since we consider exponential growth, we take the nth root of $a$:
    \[
    (a)^{\frac{1}{n}} = (\frac{(1)_n}{(\frac{1}{2})_n})^\frac{1}{n} \xrightarrow{n\rightarrow \infty} 1.
    \]
\end{proof}

\begin{Lemma}
    \label{transform}
    The following transformation holds when $b$ is non-negative:
    \[
    2F1(a,a-b+\frac{1}{2},b+\frac{1}{2},z^2)= (1+z)^{(-2a)}2F1(a,b,2b;\frac{4z}{(1+z)^2}).
    \]
\end{Lemma}
\begin{proof}
    This is a reference of \cite[Page 130, Equation (2)]{Temme-SpecialFunctions}.
\end{proof}

\begin{Proposition}
    \label{hypergeo}
    The representation gap of the Motzkin monoid has an exponential growth rate, with $n^{th}$ root $3$. (This is true over any field.)
\end{Proposition}

\begin{proof}
    The asymptotic of the formula in \autoref{simple formula} is given by Mathematica to be the Hypergeometric function:
    \[
    2F1(\frac{1}2{(-n+\sqrt{n})},\frac{1}2{(1- n+\sqrt{n})},1;4).
    \]
    From \autoref{c1/2}, we have the following:
    \[
    (\cdot 2F1(\frac{1}2{(-n+\sqrt{n})},\frac{1}2{(1- n+\sqrt{n})},1;4))^{\frac{1}{n}}
    = (2F1(\frac{1}2{(-n+\sqrt{n})},\frac{1}2{(1- n+\sqrt{n})},\frac{1}{2};4))^{\frac{1}{n}}.
    \]
    We can then apply \autoref{transform} for $a = -\frac{n}{2}+\frac{\sqrt{n}}{2}, b=0,z=2$, such that;
    \[
    (2F1(\frac{1}2{(-n+\sqrt{n})},\frac{1}2{(1- n+\sqrt{n})},\frac{1}{2};4))^{\frac{1}{n}}
    = (3^{(n-\sqrt{n})}2F1(-\frac{n}{2}+\frac{\sqrt{n}}{2},0,0,\frac{8}{9})
    )^{\frac{1}{n}}.
    \]
    From mathematica, $2F1(-\frac{n}{2}+\frac{\sqrt{n}}{2},0,0,\frac{8}{9})^{\frac{1}{n}} = 1$, where the function is defined at $c=0$ as long as $b=c \rightarrow 0$ by analytic continuation, so,
    \[
    (3^{(n-\sqrt{n})}2F1(-\frac{n}{2}+\frac{\sqrt{n}}{2},0,0,\frac{8}{9})
    )^{\frac{1}{n}} = (3^{(n-\sqrt{n})})^\frac{1}{n}.
    \]
    Now, it is easy to see that,
    \[
    \frac{3}{3^{\frac{1}{\sqrt{n}}}} \xrightarrow{n\rightarrow\infty} 3.
    \]
    This proves that the $n^{th}$ root of the representation gap converges to $3$.
\end{proof}

\subsection{Comparison of Simple and Semisimple Representations}

We now discuss where the simple and semisimple dimensions coincide, with the following data obtained from the program GAP, the code for which we reference in \autoref{gap}.

\vspace{0.2cm}

\begin{figure}[hbt!]
    \centering
    \includegraphics[scale = 0.7]{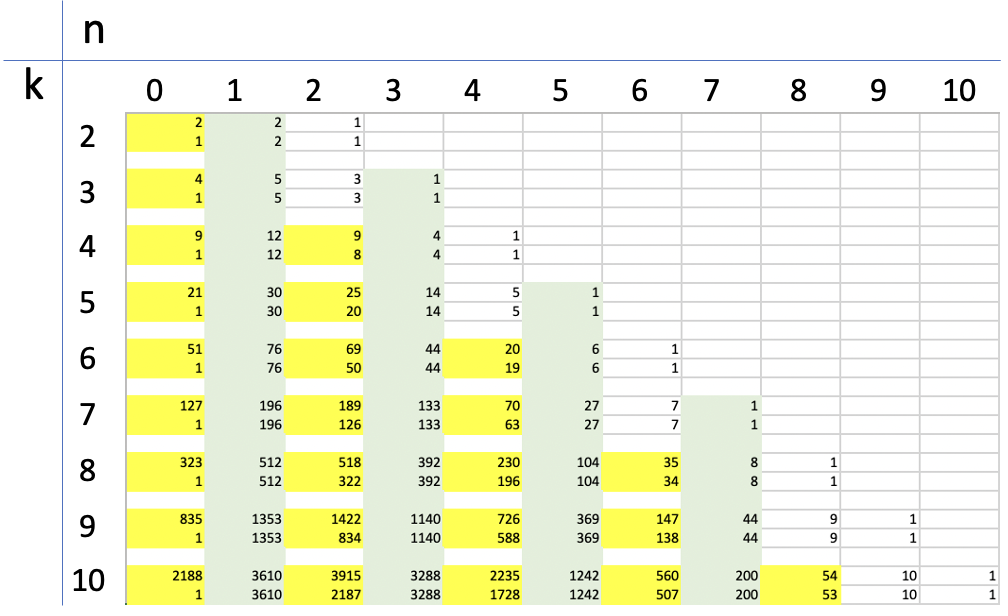}
    \caption{Semisimple and simple dimensions of the Motzkin monoid over $n$ and $k$.}
    \label{fig:faith}
\end{figure}

\autoref{fig:faith} tracks the size and rank of idempotent matrices over $n$ and $k$. The top halves of each row represent the semisimple dimension where the bottom halves state the simple dimension. Recall from \autoref{rank} that the dimension of the simple representation is the rank of the Gram matrix associated to $\jcell_k$, while the dimension of the semisimple representation is the size of the Gram matrix. We highlight columns where these values differ in yellow, and highlight in green the columns where these dimensions coincide.

We observe here that for columns with odd $k$, the Gram matrices have full rank. Indeed:

\begin{Lemma}
\label{faithful embedding}
Let $\K$ be of characteristic zero.
    For odd $K = 2k+1 \in \K$, $\dim(L_K) = \operatorname{ssdim}(L_K)$.
\end{Lemma}

\begin{proof}
    By \cite{BeHa-motzkin}, the Motzkin monoid arises as the endomorphism monoid of a tensor product of $V=\K^2\oplus\K$ (with $\K^2$ denoting the defining representation, and $\K$ the trivial representation) of quantum $\mathfrak{sl}_2$ with quantum parameter $q=\pm i$ (the imaginary unit). Now, standard results as in \cite{AnStTu-cellular-tilting} imply that the statement is equivalent to the respective highest weight summand in $V^{\otimes n}$ being simple.
    This latter fact can be deduced from the combinatorics in \cite{SuTuWeZh-mixed-tilting}.
\end{proof}

There is also a combinatorial proof in generalisation of the following example:

\begin{Example}
    \label{rowreductions}
    We exemplify this in $\mathcal{M}o_4$ for $\jcell_2$ i.e. the J-cell with two through strands. Below is a diagrammatic realisation of this $\jcell$, with idempotent elements again drawn in red:

    \begin{gather*}
\begin{tabular}{C!{\color{blue}\vrule width 0.1mm}C!{\color{blue}\vrule width 0.1mm}C!{\color{blue}\vrule width 0.1mm}C!{\color{blue}\vrule width 0.1mm}C!{\color{blue}\vrule width 0.1mm}C!{\color{blue}\vrule width 0.1mm}C!{\color{blue}\vrule width 0.1mm}C!{\color{blue}\vrule width 0.1mm}C}
\arrayrulecolor{blue}
\begin{tikzpicture}[anchorbase]
\draw[red, thick] (0,1) to (0,0);
\draw[red, thick] (0.25,1) to (0.25,0);
\draw[red, thick] (0.5,1) to[out=270,in=180] (0.625,0.75) to[out=0,in=270] (0.75,1);
\draw[red,thick] (0.5,0) to[out=90,in=180] (0.625,0.25) to[out=0,in=90] (0.75,0);
\end{tikzpicture} & 
\begin{tikzpicture}[anchorbase]
\draw[red, thick] (0,1) to (0,0);
\draw[red, thick] (0.25,1) to (0.25,0);
\draw[red, thick] (0.5,1) to[out=270,in=180] (0.625,0.75) to[out=0,in=270] (0.75,1);
\draw[red, thick] (0.5,0) to (0.5,0.25);
\draw[red, thick] (0.75,0) to (0.75,0.25);
\end{tikzpicture} &
\begin{tikzpicture}[anchorbase]
\draw[thick] (0,1) to (0,0);
\draw[thick] (0.25,1) to (0.5,0);
\draw[thick] (0.5,1) to[out=270,in=180] (0.625,0.75) to[out=0,in=270] (0.75,1);
\draw[thick] (0.25,0) to (0.25,0.25);
\draw[thick] (0.75,0) to (0.75,0.25);
\end{tikzpicture} &
\begin{tikzpicture}[anchorbase]
\draw[red, thick] (0,1) to (0,0);
\draw[red, thick] (0.25,1) to (0.75,0);
\draw[red, thick] (0.5,1) to[out=270,in=180] (0.625,0.75) to[out=0,in=270] (0.75,1);
\draw[red,thick] (0.25,0) to[out=90,in=180] (0.375,0.25) to[out=0,in=90] (0.5,0);
\end{tikzpicture} &
\begin{tikzpicture}[anchorbase]
\draw[thick] (0,1) to (0,0);
\draw[thick] (0.25,1) to (0.75,0);
\draw[thick] (0.5,1) to[out=270,in=180] (0.625,0.75) to[out=0,in=270] (0.75,1);
\draw[thick] (0.25,0) to (0.25,0.25);
\draw[thick] (0.5,0) to (0.5,0.25);
\end{tikzpicture} &
\begin{tikzpicture}[anchorbase]
\draw[thick] (0,1) to (0.25,0);
\draw[thick] (0.25,1) to (0.5,0);
\draw[thick] (0.5,1) to[out=270,in=180] (0.625,0.75) to[out=0,in=270] (0.75,1);
\draw[thick] (0,0) to (0,0.25);
\draw[thick] (0.75,0) to (0.75,0.25);
\end{tikzpicture} &
\begin{tikzpicture}[anchorbase]
\draw[thick] (0,1) to (0.5,0);
\draw[thick] (0.25,1) to (0.75,0);
\draw[thick] (0.5,1) to[out=270,in=180] (0.625,0.75) to[out=0,in=270] (0.75,1);
\draw[thick] (0,0) to[out=90,in=180] (0.125,0.25) to[out=0,in=90] (0.25,0);
\end{tikzpicture} &
\begin{tikzpicture}[anchorbase]
\draw[thick] (0,1) to (0.5,0);
\draw[thick] (0.25,1) to (0.75,0);
\draw[thick] (0.5,1) to[out=270,in=180] (0.625,0.75) to[out=0,in=270] (0.75,1);
\draw[thick] (0,0) to (0,0.25);
\draw[thick] (0.25,0) to (0.25,0.25);
\end{tikzpicture} &
\begin{tikzpicture}[anchorbase]
\draw[thick] (0,1) to (0.25,0);
\draw[thick] (0.25,1) to (0.75,0);
\draw[thick] (0.5,1) to[out=270,in=180] (0.625,0.75) to[out=0,in=270] (0.75,1);
\draw[thick] (0,0) to (0,0.25);
\draw[thick] (0.5,0) to (0.5,0.25);
\end{tikzpicture}
\\
\end{tabular}
\end{gather*}

    This has the following Gram matrix;

    \begin{center}
    \[%
,
    \]
    \end{center}

    This is not full-rank, having rank 8, but 9 rows/columns. We colour the non-linearly independent rows in blue:
\begin{gather*}
    \begin{tabular}{C!{\color{blue}\vrule width 0.1mm}C!{\color{blue}\vrule width 0.1mm}C!{\color{blue}\vrule width 0.1mm}C!{\color{blue}\vrule width 0.1mm}C!{\color{blue}\vrule width 0.1mm}C!{\color{blue}\vrule width 0.1mm}C!{\color{blue}\vrule width 0.1mm}C!{\color{blue}\vrule width 0.1mm}C}
\arrayrulecolor{blue}
\cellcolor{mydarkblue!25}
\begin{tikzpicture}[anchorbase]
\draw[red, thick] (0,1) to (0,0);
\draw[red, thick] (0.25,1) to (0.25,0);
\draw[red, thick] (0.5,1) to[out=270,in=180] (0.625,0.75) to[out=0,in=270] (0.75,1);
\draw[red,thick] (0.5,0) to[out=90,in=180] (0.625,0.25) to[out=0,in=90] (0.75,0);
\end{tikzpicture} & 
\cellcolor{mydarkblue!25}
\begin{tikzpicture}[anchorbase]
\draw[red, thick] (0,1) to (0,0);
\draw[red, thick] (0.25,1) to (0.25,0);
\draw[red, thick] (0.5,1) to[out=270,in=180] (0.625,0.75) to[out=0,in=270] (0.75,1);
\draw[red, thick] (0.5,0) to (0.5,0.25);
\draw[red, thick] (0.75,0) to (0.75,0.25);
\end{tikzpicture} &
\cellcolor{mydarkblue!25}
\begin{tikzpicture}[anchorbase]
\draw[thick] (0,1) to (0,0);
\draw[thick] (0.25,1) to (0.5,0);
\draw[thick] (0.5,1) to[out=270,in=180] (0.625,0.75) to[out=0,in=270] (0.75,1);
\draw[thick] (0.25,0) to (0.25,0.25);
\draw[thick] (0.75,0) to (0.75,0.25);
\end{tikzpicture} &
\cellcolor{pink}
\begin{tikzpicture}[anchorbase]
\draw[red, thick] (0,1) to (0,0);
\draw[red, thick] (0.25,1) to (0.75,0);
\draw[red, thick] (0.5,1) to[out=270,in=180] (0.625,0.75) to[out=0,in=270] (0.75,1);
\draw[red,thick] (0.25,0) to[out=90,in=180] (0.375,0.25) to[out=0,in=90] (0.5,0);
\end{tikzpicture} &
\cellcolor{mydarkblue!25}
\begin{tikzpicture}[anchorbase]
\draw[thick] (0,1) to (0,0);
\draw[thick] (0.25,1) to (0.75,0);
\draw[thick] (0.5,1) to[out=270,in=180] (0.625,0.75) to[out=0,in=270] (0.75,1);
\draw[thick] (0.25,0) to (0.25,0.25);
\draw[thick] (0.5,0) to (0.5,0.25);
\end{tikzpicture} &
\cellcolor{mydarkblue!25}
\begin{tikzpicture}[anchorbase]
\draw[thick] (0,1) to (0.25,0);
\draw[thick] (0.25,1) to (0.5,0);
\draw[thick] (0.5,1) to[out=270,in=180] (0.625,0.75) to[out=0,in=270] (0.75,1);
\draw[thick] (0,0) to (0,0.25);
\draw[thick] (0.75,0) to (0.75,0.25);
\end{tikzpicture} &
\cellcolor{mydarkblue!25}
\begin{tikzpicture}[anchorbase]
\draw[thick] (0,1) to (0.5,0);
\draw[thick] (0.25,1) to (0.75,0);
\draw[thick] (0.5,1) to[out=270,in=180] (0.625,0.75) to[out=0,in=270] (0.75,1);
\draw[thick] (0,0) to[out=90,in=180] (0.125,0.25) to[out=0,in=90] (0.25,0);
\end{tikzpicture} &
\cellcolor{mydarkblue!25}
\begin{tikzpicture}[anchorbase]
\draw[thick] (0,1) to (0.5,0);
\draw[thick] (0.25,1) to (0.75,0);
\draw[thick] (0.5,1) to[out=270,in=180] (0.625,0.75) to[out=0,in=270] (0.75,1);
\draw[thick] (0,0) to (0,0.25);
\draw[thick] (0.25,0) to (0.25,0.25);
\end{tikzpicture} &
\cellcolor{mydarkblue!25}
\begin{tikzpicture}[anchorbase]
\draw[thick] (0,1) to (0.25,0);
\draw[thick] (0.25,1) to (0.75,0);
\draw[thick] (0.5,1) to[out=270,in=180] (0.625,0.75) to[out=0,in=270] (0.75,1);
\draw[thick] (0,0) to (0,0.25);
\draw[thick] (0.5,0) to (0.5,0.25);
\end{tikzpicture}
\\
\end{tabular}
\end{gather*}

    Observing the corresponding diagrams, both blocks of idempotents (the 2x2 red sub-matrices in the highlighted columns) come from straight line-idempotents, which exist in the odd cases aswell, as we note that these will exist in all submonoids as described in \autoref{submonoid}. However, the two curved-idempotent diagrams (highlighted pink) which resulted in row-reductions are more interesting. These are the following:

    \begin{gather*}
    %
    \end{gather*}

    Note that they are both not only symmetric in their top row presentations, but are in fact reflections of one another about a central vertical line. This allows for the same top cell presentation to have non-unique bottom cell idempotent presentations, resulting in the row reduction (with enumerations counting rows from top to bottom): 
    \[
    \rcell_1-\rcell_2-\rcell_8+\rcell_9 = 0.
    \] 
    The general case is more difficult and omitted.
\end{Example}

In general, for $\jcell$'s with $k$ even, we always found similar patterns. Namely, we saw that an even number of strands can result in row reductions because symmetric strand positioning on the top presentation results in linear dependence in bottom presentations. As seen in \autoref{rowreductions}, this is because symmetry allows for the possibility to swap which strand is curved and which strand is straight, forming two idempotents in the same $\rcell$, resulting in possible row reductions in the Gram matrix. We hence hypothesise that $\jcell$'s with $k$ odd are linearly independent because with an odd number of through strands, such symmetric cases do not exist. However we omit a diagrammatic proof of this, and instead use the representation theoretic proof of \autoref{faithful embedding} above to generalise our observed pattern.

\subsection{Faithfulness Gap}
\label{faithful}
We now prove a lower bound for $\operatorname{faith}_\K (\mathcal{M}o_n^{\leq \sqrt{n}})$.  
Recall that faith is defined as the minimum dimension of a faithful representation of $\mathcal{M}o_n^{\leq \sqrt{n}}$. 

\begin{Proposition}
    The faithfulness has the following formula as a lower bound;
    \[
    \operatorname{faith}_\K(\mathcal{M}o_n^k) \geq \sum_{t=0}^{n}\frac{k+1}{k+t+1}\binom{n}{k+2t}\binom{k+2t}{t}, \hspace{0.1cm} \text{for k odd}.
    \]
    With the following asymptotic for odd $n$:
    \[
    \operatorname{faith}_\K(\mathcal{M}o_n^{\leq \sqrt{n}}) \geq \operatorname{ssgap}_\K(\mathcal{M}o_n^{\leq \sqrt{n}}) \geq \Theta (n^{-\frac{3}{2}}\cdot 3^n).
    \]
\end{Proposition}

\begin{proof}
    We know that the faithfulness gap is the minimum dimension of a faithful representation, and from \autoref{faithful embedding} above, we know that for odd $k$, faithful representations embed in semisimple representations, hence allowing us to use the formula for semisimple representations for odd $k$, with the growth rate we obtained in \autoref{semisimple}.
\end{proof}

\subsection{Main Results}

This brings us to our final statement, which summarises the findings of this paper:

\begin{Theorem}\label{T:Main}
    We have the following lower bounds for the exponential growth of the simple dimension, semisimple dimension and faithfulness of the Motzkin monoid; all of these are over an arbitrary field.
    \begin{enumerate}
        \item $\operatorname{gap}_\K(\mathcal{M}o_n^{\leq \sqrt{n}}) \geq f(n)3^n$ for $n\gg 0$,
        \item $\operatorname{ssgap}(\mathcal{M}o_n^{\leq \sqrt{n}}) \geq \Theta (n^{-\frac{3}{2}}\cdot 3^n)$ for $n\gg 0$,
        \item $\operatorname{faith}_\K(\mathcal{M}o_n^{\leq \sqrt{n}}) \geq \Theta (n^{-\frac{3}{2}}\cdot 3^n)$ for $n\gg 0$.
    \end{enumerate}
    Where $f(n)$ is some sub-exponential function, which we lose when taking the nth root in our calculations.
\end{Theorem}

\begin{proof}
    See \autoref{semisimple}, \autoref{simple}, and \autoref{faithful}.
\end{proof}

We now determine the ratios for each of these gaps.

\begin{Theorem}
    The Motzkin monoid has the following ratios:
    \begin{enumerate}
        \item $\lim_{n\rightarrow \infty}{(\operatorname{gapr}_{\K}(\mathcal{M}o_n^{\leq \sqrt{n}}))^{\frac{1}{n}}} = 1$
        \item $\operatorname{ssgapr}_{\K}(\mathcal{M}o_n^{\leq \sqrt{n}})\in\Theta(1)$
        \item $\operatorname{faithr}_{\K}(\mathcal{M}o_n^{\leq \sqrt{n}})\in\Theta(1)$
    \end{enumerate}
\end{Theorem}
\begin{proof}
     We use the asymptotic for the Motzkin numbers $|M| = \sqrt{\frac{3\pi}{4}}3^{n+1}n^{-\frac{3}{2}}$, which can be found in \cite{Cl-oeis-A001006}, as calculated by Benoit Cloitre. We know that the Motzkin numbers squared give the amount of Motzkin elements, so we have $|\mathcal{M}o_n|= (\sqrt{\frac{3\pi}{4}}3^{n+1}n^{-\frac{3}{2}})^2$, and we can easily obtain our ratios. Note that we ignore linear terms, since we are working with asymptotic behaviour. Consequently, using \autoref{ssdim}, we have:
    \[
    \operatorname{ssgapr}_{\K}(\monoid) = \frac{3^n n^{-\frac{3}{2}}}{\sqrt{(3^n n^{-\frac{3}{2}}})^2} = 1.
    \]
    The others, $\operatorname{gapr}_{\K}(\monoid)$ and $\operatorname{faithr}_{\K}(\monoid)$, also follow directly from \autoref{T:Main} and \autoref{ratios}.
\end{proof}
As such, we show that the representation gap ratio of the Motzkin monoid remains constant as $n$ grows very large, which is a necessary but not sufficient condition of its viability for encryption purposes.

\newcommand{\etalchar}[1]{$^{#1}$}

\appendix

\subsection{Figure 1}
Below is the mathematica code for \autoref{strandsdiagram}. This counts the number of Motzkin elements over varying parameters $n$ and $k$, then plotting a fixed $n=20$, and showing the values elements in the sum for varying values of $k \in [1,20]$.
\vspace{0.2cm}
\begin{lstlisting}[language=Mathematica]
    Mo[n_, k_] := Sum[(k + 1)/(k + t + 1) * Binomial[n, k + 2*t] * Binomial[k + 2*t, t], {t, 0, n}];
    Plot[Mo[20, k], {k, 0, 20}]
    LogPlot[Mo[20, k], {k, 0, 20}]
\end{lstlisting}

\subsection{Figure 2}
The mathematica code for \autoref{t600} again expresses the contribution of varying $t \in [0,500]$ to a sum over $t$, then plotting this.
\vspace{0.2cm}
\begin{lstlisting}[language=Mathematica]
    Moo3[n_] := 
    Table[Binomial[n, Sqrt[n] + 2*t]*Binomial[Sqrt[n] + 2*t, t], {t, 0, n}]
    ListPlot[Moo3[600], Joined -> True, PlotRange -> {{0, 500}, All}]
    ListLogPlot[Moo3[600], Joined -> True, PlotRange -> {{0, 500}, All}]
\end{lstlisting}

\subsection{Figure 3}
The mathematica code for \autoref{t=n/3} finds the maximum value of our formula over varying $t$'s, using the max argument, then comparing this with the value of $n$ divided by $3$.
\vspace{0.2cm}
\begin{lstlisting}[language=Mathematica]
    Moo[n_, k_] := DeleteCases[Table[(k + 1)/(k + t + 1)*Binomial[n, k + 2*t]*
    Binomial[k + 2*t, t], {t, 0, n/2}], 0, Infinity];
    
    max[n_] := Module[{}, NS = Moo[n, Sqrt[n]];
    
    NN := Table[Total[NS[[b]]], {b, 1, Length[NS]}];
    
    Position[NN, Max[NN]][[1, 1]]];
    
    X = Table[max[i], {i, 1, 10000}];
    
    ListPlot[{X, Table[0.333333*k, {k, 1, 10000}]}]
    ListLogPlot[{X, Table[0.333333*k, {k, 1, 10000}]}]
\end{lstlisting}

\subsection{Semisimple Gap}
Here, we input our revised formula for the semisimple gap, then using the Asymptotic argument in Mathematica to obtain an approximation for our growth rate, applying this argument twice to simplify as much as possible, as we are mostly interested in the exponential growth.
\vspace{0.2cm}
\begin{lstlisting}[language=Mathematica]
    formula[n_, a_] := ((n^(1/2) + 1)/(n^(1/2) + 1 + n/a))* Binomial[n, n/a] *Binomial[n - n/a, n^(1/2) + n/a];
    
    Asymptotic[formula[n, 3], n -> Infinity]
    
\end{lstlisting}
\subsection{Representation Gap}
This mathematica code works similarly as above, the first part is the code for the plot which we produced, the second part gives us the hypergeometric function as our asymptotic, and the third tests values of the hypergeometric function when its second and third arguments are 0.
\vspace{0.2cm}
\begin{lstlisting}[language=Mathematica]
    Plot[{Hypergeometric2F1[1/2*(-y + Sqrt[y]), 1/2*(1 - y + Sqrt[y]), 1, 4]^(1/y), 3, 2}, {y, 10, 5000},WorkingPrecision -> 30]
    
    Asymptotic[Sum[Binomial[i - Floor[Sqrt[i]], 2 t] Binomial[2 t, t], {t, 1, i}], i -> Infinity] // Simplify
    
    Hypergeometric2F1[a, 0, 0, b]
\end{lstlisting}

\subsection{Simple and Semisimple Dimensions}
\label{gap}
The following code was input into Gap in order to obtain values for our simple and semisimple dimensions. The first part sets up our functions and the second outputs dimensions and ranks of Gram matrices over varying $k$, inducting on $n$. Note that due to exponential growth in size, obtaining these results takes a longer time for each $n+1$, making values higher than roughly $n=12$ difficult to obtain.
\vspace{0.2cm}
\begin{lstlisting}[language=GAP]
01Matrix := function(D)
  local result, row, i;
  result := List(Matrix(PrincipalFactor(D)), ShallowCopy);
  for row in result do
    for i in [1 .. Length(row)] do
      if row[i] <> 0 then
        row[i] := 1;
      fi;
    od;
  od;
  return result;
end;

DClassMatrixRanks := function(S)
  local D, n, L;
  D := DClasses(S);
  n := Length(D);
  L := List([1..n], i -> [Length(01Matrix(D[i])),
                          Rank(01Matrix(D[i]))]);
  return L;
end;

gap> n := 1; DClassMatrixRanks(MotzkinMonoid(n));
gap> n := n + 1; DClassMatrixRanks(MotzkinMonoid(n));
\end{lstlisting}

\end{document}